\documentclass[12pt]{article} 
\usepackage{setspace}
\usepackage[utf8]{inputenc} 

\usepackage{subfigure}
\usepackage[margin=1in]{geometry}
\usepackage{graphicx} 
\usepackage{epstopdf}

\usepackage{booktabs} 
\usepackage{array}
\usepackage{paralist} 
\usepackage{verbatim} 
\usepackage{subfig}

\usepackage{fancyhdr} 
\usepackage{sectsty}
\allsectionsfont{\sffamily\mdseries\upshape} 

\usepackage[nottoc,notlof,notlot]{tocbibind} 
\usepackage[titles,subfigure]{tocloft}

\usepackage{amsmath, amsthm, amssymb}

\usepackage{multirow}
\usepackage{algorithm}
\usepackage{algorithmic}

\usepackage{amsthm}
\newtheorem{thm}{Theorem}[section]

\usepackage{rotating}
\usepackage{pdflscape}
 \usepackage{url}

\begin{document}

\title{ {Inexactness of SDP Relaxation and Valid Inequalities for Optimal Power Flow}}

\author{Burak Kocuk, Santanu S. Dey, X. Andy Sun  
\thanks{
The authors are with the School of Industrial and Systems Engineering, Georgia Institute of Technology, 765 Ferst Drive, NW Atlanta, Georgia 30332-0205 (e-mail: burak.kocuk, santanu.dey, andy.sun@isye.gatech.edu).
}
}
\vspace{-5mm}

\markboth{SUBMITTED TO IEEE TRANSACTIONS ON POWER SYSTEMS, September 2014}%
{Shell \MakeLowercase{\textit{et al.}}: Bare Demo of IEEEtran.cls for Journals}
\maketitle

\begin{abstract}
It has been recently proven that the semidefinite programming (SDP) relaxation of the optimal power flow problem over radial networks is exact {under technical conditions such as not including  generation lower bounds or allowing load over-satisfaction}. In this paper, we investigate the situation where generation lower bounds are present. {We show that even for a two-bus one-generator system, the SDP relaxation can have all possible approximation outcomes, that is (1) SDP relaxation may be exact or (2) SDP relaxation may be inexact or (3) SDP relaxation may be feasible while the OPF instance may be infeasible. We provide a complete characterization
of when these three approximation outcomes occur and an analytical expression of the resulting optimality gap for this two-bus system.} {In order to facilitate further research, we design a library of instances {over radial networks} in which the SDP relaxation has positive optimality gap}. {Finally,} 
we propose valid inequalities and variable bound tightening techniques 
{that significantly improve} the computational performance of a global optimization solver. Our work demonstrates the need of developing efficient global optimization methods  for the solution of OPF even in the {simple but fundamental} case of radial networks.
\end{abstract}



\section{Introduction}

Optimal Power Flow (OPF) was first introduced in the 1960s \cite{Carpentier}  and much effort has been devoted to its solution, which {has} resulted in a rich literature. Roughly speaking, we can categorize the previous work into three categories.
 
The first {category} of algorithms find local optimal solutions or stationary points using optimization procedures such as interior point methods  (e.g. MATPOWER \cite{Matpower}). {The shortcoming of these local methods is that if a solution is obtained, there is  no guarantee for global optimality or even any evidence of how good this solution is.} {For instance, in \cite{Bukhsh}, there are several examples which have multiple local optima and it has been shown that local solvers tend to converge to the solution which is closest to the initial guess.}

The second {category} of algorithms attempt to obtain global optimal solutions of OPF by solving convex relaxations. For instance, one popular approach is to use {semidefinite programming (SDP)} relaxations of the original OPF problem. Since SDPs are polynomially solvable, this method suggests that OPF can be solved efficiently provided that the relaxation is exact, i.e., the SDP relaxation finds the global optimal solution for the OPF problem. {A recent, comprehensive survey can be found in \cite{low2014convexi, low2014convexii}.}

\begin{itemize}
\item[-]  One of the early works that popularizes {this} approach is \cite{Lavaei12}.
{It is shown that the SDP relaxation is tight for a resistive network with no reactive loads where demand over-satisfaction is allowed, as long as the dual variables are positive.}
{It was conjectured that under normal operating conditions the SDP relaxation is tight.}
\item[-] {However, \cite{Lesieutre} gives a very simple counterexample (a 3-bus cycle) with nonzero optimality gap.}

\item[-] {In \cite{sojoudi2012physics}, it is proven that SDP relaxation is exact if load over-satisfaction is allowed and a sufficient number of virtual phase shifters are present.}

  

\item[-] An attempt to solve OPF using SDP relaxation is made in \cite{Zhang12} for radial networks.  In this work, it is proven that under operational constraints on voltage magnitudes, line losses, and line flows, the SDP relaxation is tight {if there are no lower bounds on real and reactive power generation at any bus}. {Similar results are also obtained in \cite{bose2011optimal, bose2012quadratically} without line limit constraints.}

\item[-] In \cite{Lavaei14}, it is proven that if voltage magnitudes are fixed, then the convex relaxations are tight under practical angle restrictions {for radial networks} {in the presence of {only} real power lower bounds.} This result extends to the case with variable voltage magnitudes under reasonable assumptions.
\end{itemize}

As we observe above, the exactness of the SDP relaxation can only be guaranteed for special classes of OPF instances, {often when we disregard some generation lower bounds}. {Unfortunately,} if the SDP relaxation is not tight, the physical meaning of its solution is not easy to recover. {In this case, an alternative approach would be to use an hierarchy for SDPs \cite{lasserre}  as suggested in \cite{josz}.  However, this approach may not be practical as the size of the SDPs grows larger with the order of the relaxation.}

{The third category of algorithms attempt to remove the pitfalls of the previous two approaches by endeavoring to obtain globally optimal solutions. One such algorithm based on branch-and-bound method is proposed in \cite{Phan} for the solution of OPF.} Lagrangian relaxation is used to find lower bounds while a local solver (IPOPT) is utilized to obtain upper bounds. {Global solution techniques are in their infancy today and much work needs to be done to make them practically efficient.}

In this paper, we focus on {the OPF problem on }radial networks {in the presence of}
generation lower bounds {on both real and reactive power.} {The goal of this paper is two fold: To highlight the inexactness of standard convex relaxations for these instances and to make algorithmic progress in solving such instances globally.}
We make two comments here in relation to the class of OPF problems we consider and our assumptions. First, although most power flow networks are not radial, they are {usually} quite sparse and analyzing radial networks can {therefore} be beneficial for their own right, especially in the case of distribution networks \cite{Lavaei14}.
{Second, typically power systems have ramping constraints, so that the power generation in the next time period cannot deviate from the current one too much. Hence, it is important to make a study of the effects of lower bounds.}



In practice, SDPs may become prohibitively expensive as the size of the network grows larger. One can turn to {second order conic programming (SOCP)} relaxations, which are in general weaker than their SDP counterparts. However, in  \cite{Sojoudi}, it has been proven that both types of relaxations give the same lower bound for the OPF problem over radial networks even if they are inexact. {Therefore any result stated for SOCPs in this paper holds for SDP relaxations and vice-versa.}

The rest of the paper is organized as follows: In the next section, we review the rectangular formulation of the OPF problem and a reformulation that leads to the SOCP relaxation. In Section \ref{section:twobus}, we {begin working on our first goal by providing}  a complete characterization of the approximation performance of SOCP relaxation for a two-bus system. In Section \ref{section:examples}, we further study the feasible regions of two small systems. Then, in Section \ref{section:library}, {we begin working on our second goal by providing} a library of radial network instances generated from MATPOWER test cases
for which SOCP relaxation is inexact. 
In Section \ref{section:valid}, we propose valid inequalities for the SOCP relaxation, which significantly improves the computational performance of a global solver. {Concluding remarks} are {made} in Section \ref{section:conc}.


\section{Optimal Power Flow}
Consider a typical power network, where $\mathcal{B}$, $\mathcal{G}$ and $\mathcal{L}$ denote respectively the set of buses, generators and transmission lines.
The nodal admittance matrix $Y \in \mathbb{C}^{|\mathcal{B}| \times |\mathcal{B}|}$ has component $Y_{ij}=G_{ij} + \mathrm{i}B_{ij}$ for each line $(i,j)$ and $G_{ii}=g_{ii}-\sum_{j\ne i} G_{ij}, B_{ii}=b_{ii}-\sum_{j\ne i} B_{ij}$, where $g_{ii}$ (resp. $b_{ii}$) is the shunt conductance (resp. susceptance) at bus $i$. 
Let $p_i^g, q_i^g$ (resp. $p_i^d, q_i^d$) be the real and reactive power output of the generator (resp. load) at bus $i$. The complex voltage $V_i$ at bus $i$ can be expressed either in the rectangular form as $V_i = e_i+\mathrm{i} f_i$ or in the polar form as $V_i = |V_i|(\cos\theta_i+\mathrm{i}\sin\theta_i)$, where the voltage magnitude $|V_i|^2=e_i^2 + f_i^2$. 

The OPF problem in the rectangular form is given as
\begin{subequations} \label{rect form}
\begin{align}
  \min  &\  \sum_{i \in \mathcal{G}} C_i(p_i^g)  \label{objR} \\
  \mathrm{s.t.} \  & p_i^g-p_i^d =\sum_{j \in \delta(i)}[ G_{ij}(e_ie_j+f_if_j)-B_{ij}(e_if_j-e_jf_i)]   \notag\\
                          &     \hspace{61mm}\quad i \in \mathcal{B} \label{activeAtBus}\\
  &  q_i^g-q_i^d = \sum_{j \in \delta(i)}[ -B_{ij}(e_ie_j+f_if_j)-G_{ij}(e_if_j-e_jf_i)]   \nonumber  \\
  & \hspace{65.5mm} i \in \mathcal{B} \label{reactiveAtBus}\\
  & (V_i^{\text{min}})^2 \le e_i^2+f_i^2 \le (V_i^{\text{max}})^2    \quad\quad\quad\quad\quad\quad \;  i \in \mathcal{B} \label{voltageAtBus} \\
  &   p_i^{\text{min}}  \le p_i^g \le p_i^{\text{max}}    \quad\quad\quad\quad\quad\quad\quad\quad\quad\quad\quad \; \; { i \in \mathcal{B}} \label{activeAtGenerator} \\
  &  q_i^{\text{min}}  \le q_i^g \le q_i^{\text{max}}    \quad\quad\quad\quad\quad\quad\quad\quad\quad\quad\quad \ \;  {i \in \mathcal{B}} \label{reactiveAtGenerator}
\end{align}
\end{subequations}
Here $C_i(p_i^g)$ in \eqref{objR} represents the production cost of generator $i$, which typically is either a linear or a convex quadratic nondecreasing function of $p_i^g$. Constraints \eqref{activeAtBus}-\eqref{reactiveAtBus} enforce flow conservation at each bus $i$, where $\delta(i)$ is the set of buses adjacent to $i$ and including $i$. Constraint \eqref{voltageAtBus} limits the upper and lower bounds on the bus voltage magnitudes. Usually $V_i^{\text{min}}$ and $V_i^{\text{max}}$ are close to the unit voltage. Constraints \eqref{activeAtGenerator}-\eqref{reactiveAtGenerator} are the upper and lower bounds on generator $i$'s real and reactive power, respectively. {Here, we have $ p_i^{\text{min}} =  p_i^{\text{max}}=q_i^{\text{min}} =  q_i^{\text{max}}=0$ for bus $i$ where there is no generator, i.e. $i \in \mathcal{B} \setminus \mathcal{G}$.}

One can equivalently formulate the above OPF problem in polar coordinates. Sometimes, the rectangular formulation is preferred since the Hessian matrix of the constraints is constant and this is an advantage for the interior point methods. On the other hand, when the voltage magnitude is fixed at some buses, the polar formulation may become more advantageous \cite{Sun}.


We can observe that all the nonlinearities in \eqref{rect form} are of the following three types:
\begin{align*}
(1)\; e_i^2+f_i^2\quad (2)\; e_ie_j+f_if_j\quad (3)\; e_if_j-e_jf_i,
\end{align*}
which are equal to $|V_i|^2$, $|V_i||V_j|\cos(\theta_i - \theta_j)$, and $|V_i||V_j|\sin(\theta_i - \theta_j)$ in the polar form, respectively. Let us define new variables $c_{ii}$, $c_{ij}$, and $s_{ij}$ for each of these three quantities. Since the cosine function is even and the sine function is odd, we also have $c_{ij}=c_{ji}$  and $s_{ij}=-s_{ji}$. On each line $(i,j)$, these quantities are linked through the fundamental trigonometric identity $\cos^2(\theta_i-\theta_j)+\sin^2(\theta_i-\theta_j)=1$, which translates into
\begin{align*}
(e_ie_j+f_if_j)^2 + (e_if_j-e_jf_i)^2  &= (e_i^2+f_i^2)(e_j^2+f_j^2)
\end{align*}
in the rectangular form. In the space of our new variables, this relation is expressed in the following quadratic equation $c_{ij}^2+s_{ij}^2  = c_{ii}c_{jj}$,
which describes the surface of a rotated second-order cone in four dimensions.

Now, we are ready to reformulate OPF using this idea:
\begin{subequations} \label{SOCP}
\begin{align}
  \min  &\ \sum_{i \in \mathcal{G}} C_i(p_i^g)  \label{obj} \\
  \mathrm{s.t.} \  &\hspace{0.5em} p_i^g-p_i^d = 
\sum_{j \in \delta(i)}[ G_{ij}c_{ij} -B_{ij}s_{ij}]  \hspace{14mm} i \in \mathcal{B} \label{activeAtBusR} \\
 & q_i^g-q_i^d = 
\sum_{j \in \delta(i)}[ -B_{ij}c_{ij} -G_{ij}s_{ij}] \hspace{14mm} i \in \mathcal{B} \label{reactiveAtBusR} \\
  &  (V_i^{\text{min}})^2 \le c_{ii} \le (V_i^{\text{max}})^2  \hspace{2.8cm}  i\in \mathcal{B} \label{voltageAtBusR} \\
  &   p_i^{\text{min}}  \le p_i^g \le p_i^{\text{max}}  \hspace{3.9cm}    {i \in \mathcal{B}}  \label{activeAtGeneratorR} \\
  &   q_i^{\text{min}}  \le q_i^g \le q_i^{\text{max}}   \hspace{3.9cm}   { i \in \mathcal{B}} \label{reactiveAtGeneratorR} \\
  & c_{ij}=c_{ji} \hspace{4.5cm}   (i,j) \in \mathcal{L} \label{cosine}\\
  & s_{ij}=-s_{ji}\hspace{4.2cm}   (i,j) \in \mathcal{L} \label{sine}\\
  & c_{ij}^2+s_{ij}^2  = c_{ii}c_{jj}. \hspace{3.2cm}     (i,j) \in \mathcal{L} \label{coupling}
\end{align}
\end{subequations}
This reformulation \eqref{SOCP} is exact for any radial network, because the following equations on voltage angles
\begin{align}
{\sin(\theta_i - \theta_j) = \frac{s_{ij}}{\sqrt{c_{ii}c_{jj}}}, \hspace{1cm} (i,j) \in \mathcal{L} \label{thetasine}}\\
{\cos(\theta_i - \theta_j) = \frac{c_{ij}}{\sqrt{c_{ii}c_{jj}}}, \hspace{1cm} (i,j) \in \mathcal{L} \label{thetacosine}}
\end{align}
have a unique solution as long as the underlying network is radial, where $s_{ij}, c_{ij}$ are obtained from solving \eqref{SOCP}. {An alternative proof  can be seen in \cite{Gan}.} For meshed networks, however, the reformulation (\ref{SOCP}) is exact only if we include \eqref{thetasine}-\eqref{thetacosine} in the constraints. 
This idea is first proposed to solve the load flow problem for radial and meshed networks \cite{Exposito99, Jabr06,Jabr07}. Then, it is adapted to OPF in \cite{Jabr08}.


Except the coupling constraints  (\ref{coupling}), all other constraints in \eqref{SOCP} are now linear. Hence, all the nonconvexity of the OPF problem \eqref{rect form} in a radial network is captured by \eqref{coupling}, and the feasible region is the intersection of a polytope defined by  (\ref{activeAtBusR})-(\ref{sine}) with the boundaries of rotated second-order cones defined by (\ref{coupling}).
It is straightforward to obtain a second-order cone programming (SOCP) relaxation of \eqref{SOCP} by relaxing  constraint (\ref{coupling}) as follows:
\begin{equation}
c_{ij}^2+s_{ij}^2  \le  c_{ii}c_{jj}  \quad (i,j) \in \mathcal{L}, \label{couplingSOCP}
\end{equation}
which can be written more explicitly as a SOCP constraint:
\begin{equation}
c_{ij}^2+s_{ij}^2 + \left(\frac{c_{ii }-c_{jj}}{2} \right)^2 \le  \left(\frac{c_{ii }+ c_{jj}}{2} \right)^2 \quad (i,j) \in \mathcal{L}.\label{SOCPconstr}
\end{equation}
The SOCP relaxation is defined as  \eqref{obj}-\eqref{sine} and \eqref{SOCPconstr}.
It is proven that in radial networks the SOCP relaxation is equivalent to the SDP relaxation \cite{Sojoudi}. In this paper, we focus on SOCP relaxation due to its superior computational performance.

\section{Analytical Study of a Two-Bus System}\label{section:twobus}


In this section, we study the two-bus system with one generator and one load. This is arguably the simplest power system, but also one of the most fundamental models in power system analysis. Surprisingly, for this simple system, the SOCP relaxation with generation lower bounds can have all three possible outcomes in terms of optimality gap, namely (1) SOCP obtains exact solution (i.e. optimality gap is zero); (2) SOCP is feasible, yet OPF is infeasible (optimality gap is infinite); (3) SOCP has a finite optimality gap, and we give an analytical expression of this gap. We identify parameter ranges in closed form for each of these outcomes. We also study the feasible region projected in the space of squared bus voltage magnitudes to gain geometric intuition.

\begin{figure*}
\begin{center}
\caption{
Projection of feasible region of 2-bus, 1-generator examples onto $(c _{11}, c_{22})$ space for five cases.  Horizontal axis is $c _{11}$ and vertical axis is $c _{22}$. Solid black curve is {(\ref{feasCurve}) containing the} feasible region of OPF with dashed lines being two asymptotes shown in Fig. 1a. Green and red lines are bounds on $c _{11}$ and $c _{22}$, resp. Magenta line is the effective lower bound on $c _{11} - c _{22}$. Blue region is the feasible region of SOCP relaxation. All figures are in p.u.
}  \label{Cases}
	\subfigure[Case 1: SOCP relaxation is exact.]{
            \includegraphics[width=0.356\columnwidth]{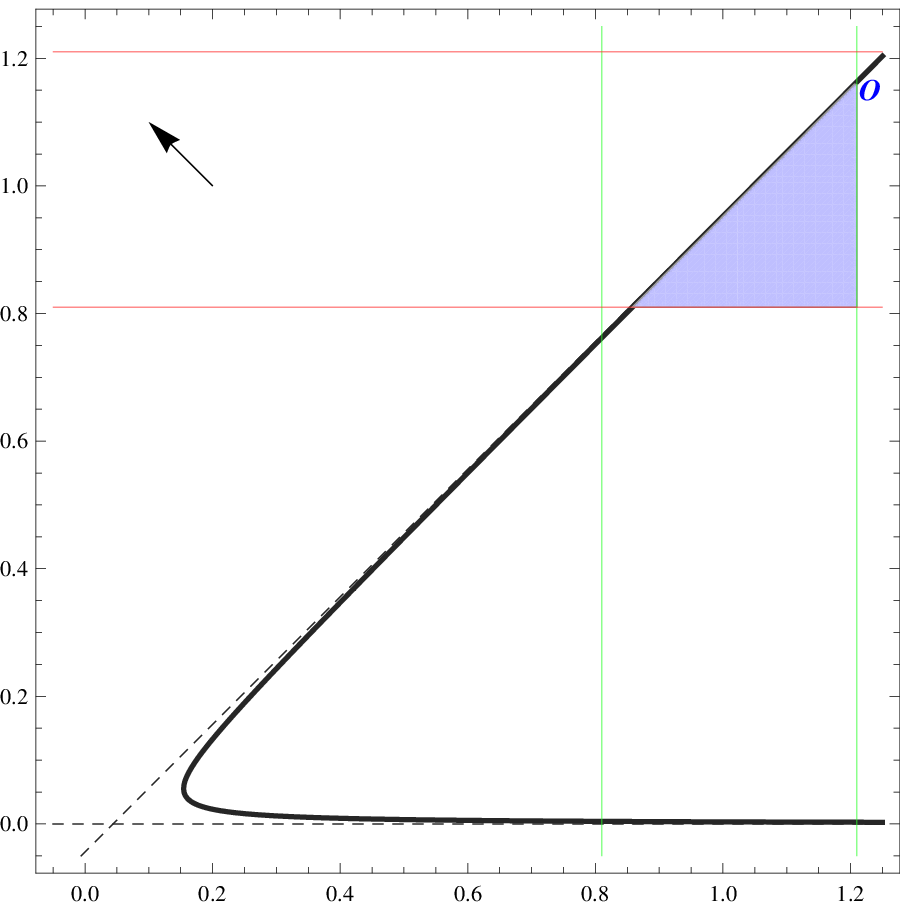}
	}
	\subfigure[Case 2: SOCP relaxation is feasible, OPF is infeasible.]{
            \includegraphics[width=0.356\columnwidth]{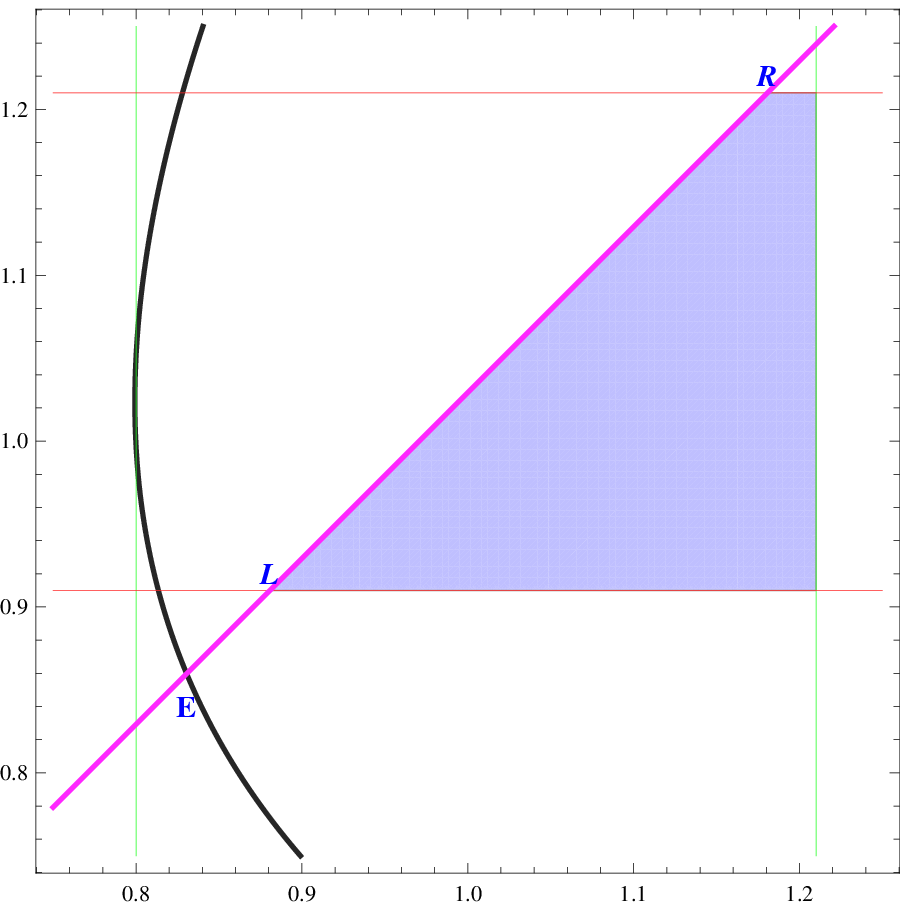}
	}
	\subfigure[Case 3: SOCP relaxation is exact.]{
            \includegraphics[width=0.356\columnwidth]{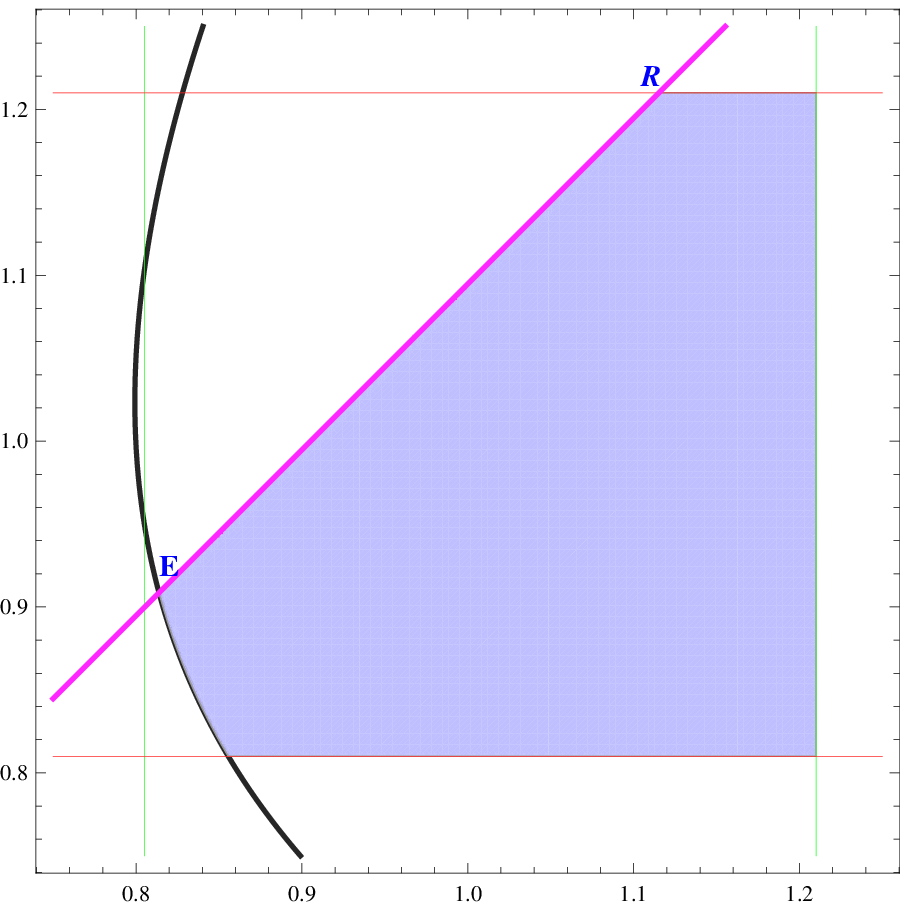}
	}
	\subfigure[Case 4: SOCP relaxation is feasible, OPF is infeasible.]{
            \includegraphics[width=0.356\columnwidth]{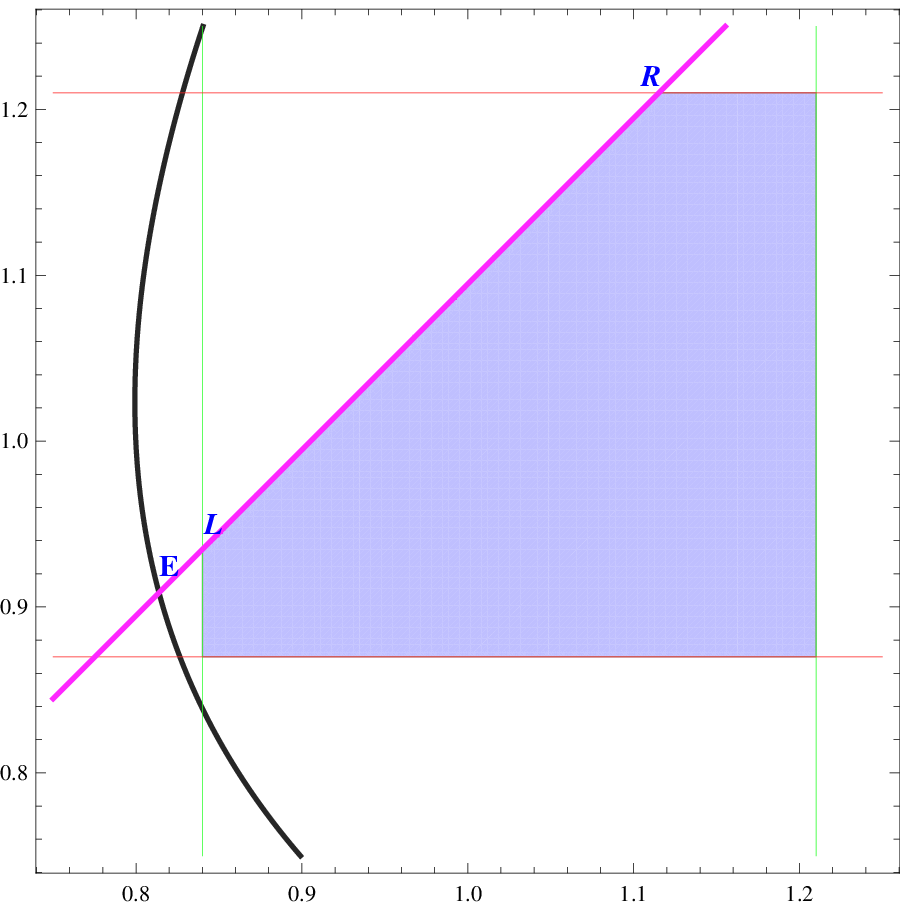}
	}
	\subfigure[Case 5: SOCP relaxation is inexact.]{
            \includegraphics[width=0.356\columnwidth]{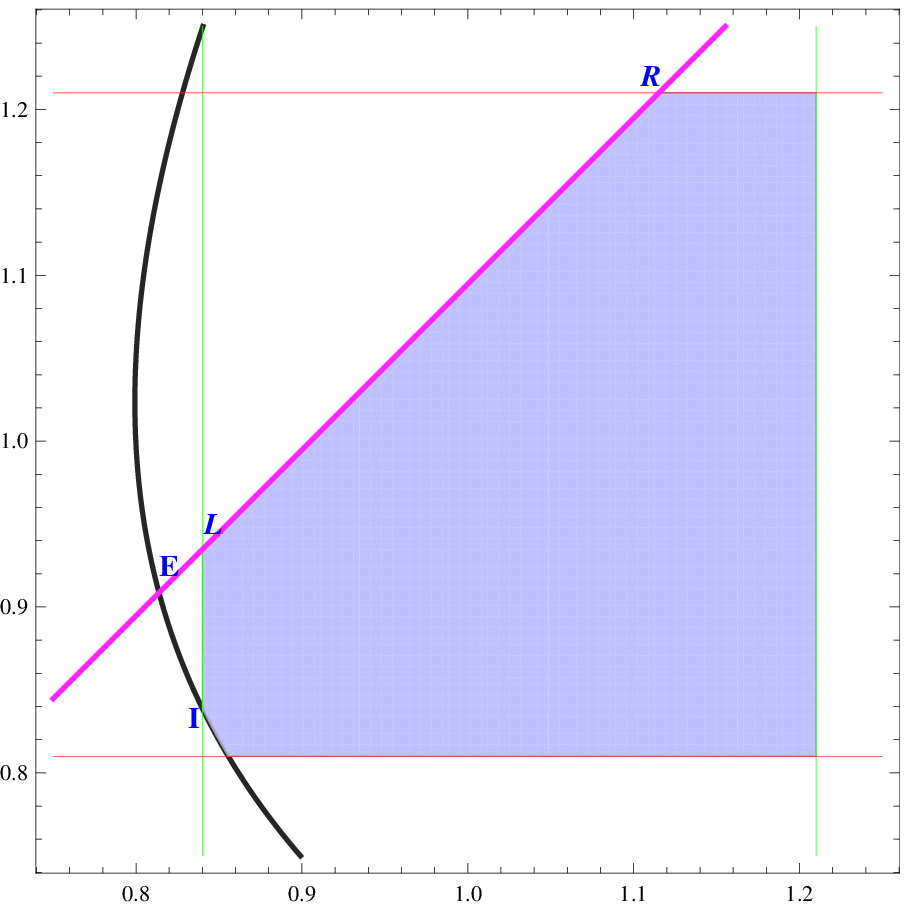}
	}
\end{center}
\end{figure*}

Let us assume that bus 1 is a generator bus and bus 2 is a  load bus. Further assume that $g_{ii}=b_{ii}=0$ and $G:= G_{12}<0$ and $B:= B_{12}>0$ (the analysis for $B<0$ is similar). Also assume the production cost $C_1(p_1^g)$ is linear in $p_1^g$.

\subsection{Feasible Region Projected to $(c_{11},c_{22})$ space}
The linear equality system (\ref{activeAtBusR})-(\ref{reactiveAtBusR}) can be written as
\begin{align}
\begin{bmatrix}
 1 &  &  \ \ G &    & -G & \ \ B   \\
 &  1  & -B  &  &\ \ B &  \ \ G  \\
& &    & \ \  G & -G& -B  \\
&  &    & -B &\ \  B&    -G
\end{bmatrix}
\begin{bmatrix}
p_1^g   \\
q_1^g   \\
c_{11}  \\
c_{22}  \\
c_{12}  \\
s_{12}
\end{bmatrix}
=
\begin{bmatrix}
0  \\
0   \\
p_2^d   \\
q_2^d
\end{bmatrix}.\label{eq:linsys}
\end{align}
Let us define
\begin{equation}
\alpha := \frac{B p_2^d + G q_2^d}{B^2+G^2}\;\;
\text{ and }\;\;
\beta :=  \frac{G p_2^d - B q_2^d}{B^2+G^2},
\end{equation}
which are constant for fixed $B,G$ and load. Solving the linear system \eqref{eq:linsys}, we can express $(p_1^g, q_1^g, c_{12}, s_{12})$ in terms of $(c_{11}, c_{22})$ as follows
\begin{subequations}
\begin{align}
s_{12} &=  -\alpha \label{eq:s12}\\
c_{12} &= c_{22}  - \beta \label{eq:c12}\\
p_1^g  &= -G(c_{11} - c_{22})  - G\beta + B \alpha \label{eq:p1g}\\
q_1^g  &= B(c_{11} - c_{22})   + B\beta + G \alpha. \label{eq:q1g}
\end{align}
\end{subequations}
We {now} reformulate constraint (\ref{coupling}) {using (\ref{eq:s12}) and (\ref{eq:c12})} as
\begin{equation} \label{feasCurve}
(c_{22} - \beta)^2 + \alpha^2 = c_{11} c_{22} \Rightarrow c_{11} = c_{22} - 2\beta + \frac{\alpha^2+\beta^2}{c_{22}},
\end{equation}
which defines a hyperbola for $(c_{11}, c_{22})$ with two asymptotes:  $c_{11} - c_{22}= - 2\beta$ and $c_{22} = 0$.

Observe that this hyperbola together with the constraints on $c_{11}$ and $c_{22}$ implied from \eqref{voltageAtBusR}-\eqref{reactiveAtGeneratorR} define the feasible region of the OPF problem projected to the $(c_{11}, c_{22})$ space. In particular, \eqref{voltageAtBusR} impose a box constraint on $c_{11}$ and $c_{22}$, whereas \eqref{activeAtGeneratorR}-\eqref{reactiveAtGeneratorR} imply upper and lower bounds on the difference $c_{11}-c_{22}$, which defines a region parallel to the first asymptote $c_{11}-c_{22}=-2\beta$. {Figure 1a depicts the entire feasible regions of OPF in black curve and of SOCP relaxation in the blue region. Figures 1b-1e zoom in particular parts. }

Furthermore, since the objective function $C_1(p_1^g)$ is assumed to be linear in $p_1^g$ and by \eqref{eq:p1g}, we can see that the level set of the objective function in $(c_{11}, c_{22})$ is also parallel to the first asymptote, and decreases toward the upper left corner as pointed by the arrow in Figure 1. Therefore, only the lower bounds on $p_1^g$ and $q_1^g$ can affect the optimal solution of \eqref{SOCP}. For this reason, we find the \textit{effective} lower bound for the difference $c_{11}-c_{22}$ as
\begin{equation}
\Delta =  \max\left\{ \frac{ p_1^{\text{min}} + G\beta - B \alpha}{-G}, \frac{q_1^{\text{min}}  - B\beta - G \alpha }{B}\right\},\label{eq:Delta}
\end{equation}
which is given by the lower bounds of \eqref{activeAtGeneratorR}-\eqref{reactiveAtGeneratorR}, and is plotted as magenta lines in Figure 1. Also note that as  $p_1^{\min}$ and $q_1^{\min}$ increase, the line $c_{11}-c_{22}\geq \Delta$ moves toward the lower right corner in Figure 1.


\subsection{Complete Characterization of Approximation Outcomes}
At this point, we are ready to explore the optimal solutions of the OPF \eqref{SOCP} and its SOCP relaxation and classify all five possible cases of the configurations of their feasible regions and the associated approximation outcomes.

\begin{itemize}
\item First of all, let us assume that $\Delta$ defined in \eqref{eq:Delta} is small enough. In this case, as depicted in Figure \ref{Cases}a, the optimal solution of both the OPF and the SOCP is unique and given by
\begin{align}
\hspace{0mm}(c_{11}^O, c_{22}^O) = \begin{cases}
\bigl(\overline c_{22} - 2\beta + \frac{\alpha^2+\beta^2}{ \overline c_{22}}, \overline c_{22}\bigr) \quad \text{ if (a) holds}   \\
\bigl(\overline c_{11}, \frac{2\beta + \overline c_{11} + \sqrt{(2\beta + \overline c_{11})^2 - 4(\alpha^2+\beta^2)}}{2}\bigr)\text{ o.w.}
\end{cases}\label{eq:cO}
\end{align}
where condition (a) is  $(\overline c_{22} - \beta)^2 + \alpha^2 \le \overline c_{11} \overline c_{22}$ and $\overline c_{ii}:=(V_i^{\max})^2$.
Hence, the SOCP relaxation is exact. This result is in accordance with the results in \cite{Zhang12}.

\item Consider the case where  $\Delta$ is large enough. In particular, $c_{11}^O - c_{22}^O < \Delta$.
Define the intersection of $c_{11} - c_{22} = \Delta$ with the binding upper bound of either $c_{11}$ or $c_{22}$ as
\begin{equation}
(c_{11}^R, c_{22}^R) = \begin{cases}   ( \overline c_{22} + \Delta, \overline c_{22} ) \text{ if }   \overline c_{11} -\overline c_{22} \ge \Delta  \\
(\overline c_{11}, \overline c_{11} - \Delta) \text{ o.w. }
\end{cases}
\end{equation}
Note that this point is OPF infeasible despite being  SOCP optimal. Next, define the intersection of the hyperbola (\ref{feasCurve}) and $c_{11} - c_{22} = \Delta$ as
\begin{equation}
(c_{11}^E, c_{22}^E) = \left( \frac{\alpha^2+\beta^2}{2\beta+  \Delta } + \Delta, \frac{\alpha^2+\beta^2}{2\beta+  \Delta }\right).\label{eq:cE}
\end{equation}

\begin{figure*}
\begin{center}
\caption{Projection of feasible region of 2-bus, 2-generator example onto $(p_1^g,p_2^g)$ space. Horizontal axis is $p_1^g$ and vertical axis is $p_2^g$.
Black curve is an ellipse with counterclockwise orientation that {contains} the feasible region of OPF problem whereas blue region is its SOCP relaxation.  Green lines are the lower bound on $p_1^g$ and $p_2^g$ while red lines are the lower bound on $q_1^g$ and $q_2^g$. Dashed lines represent angle bounds corresponding to $30^{\circ}$. Assuming linear functions, the arrow shows the cost vector. Blue and orange dots are respectively the optimal solutions of SOCP relaxation and OPF, whenever the latter exists.}  \label{FeasReg2 all}
	\subfigure[SOCP is exact.]{
            \includegraphics[width=0.465\columnwidth]{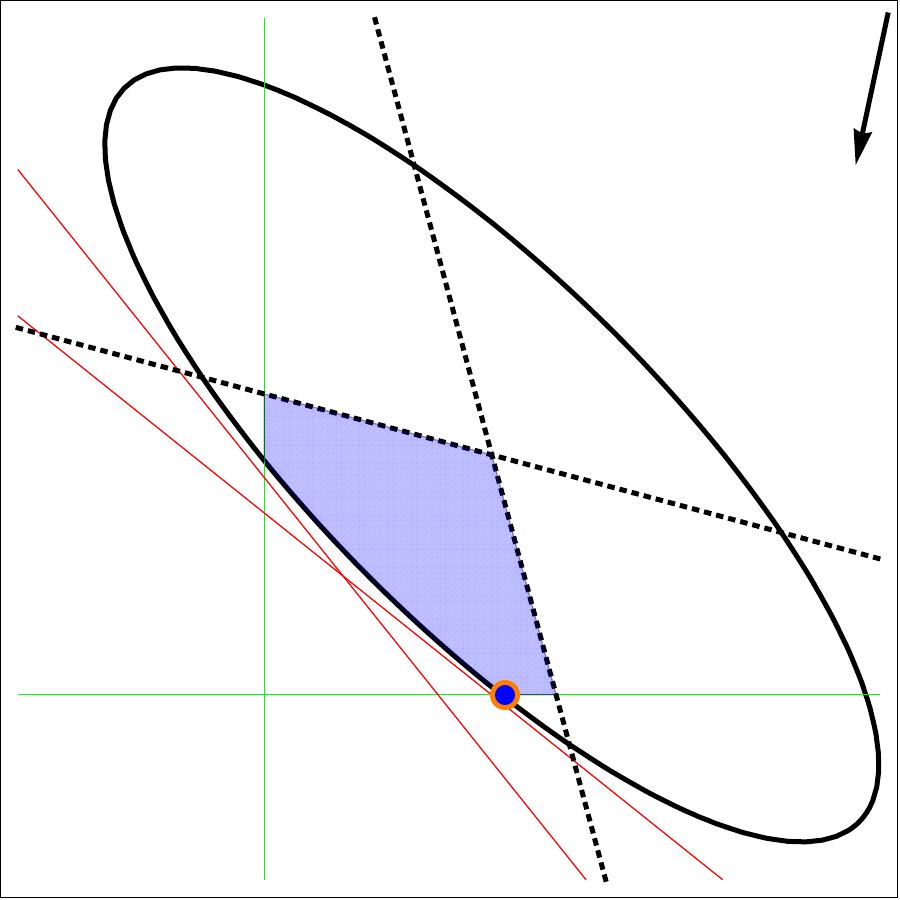}
	}
	\subfigure[SOCP is feasible while OPF is infeasible.]{
            \includegraphics[width=0.465\columnwidth]{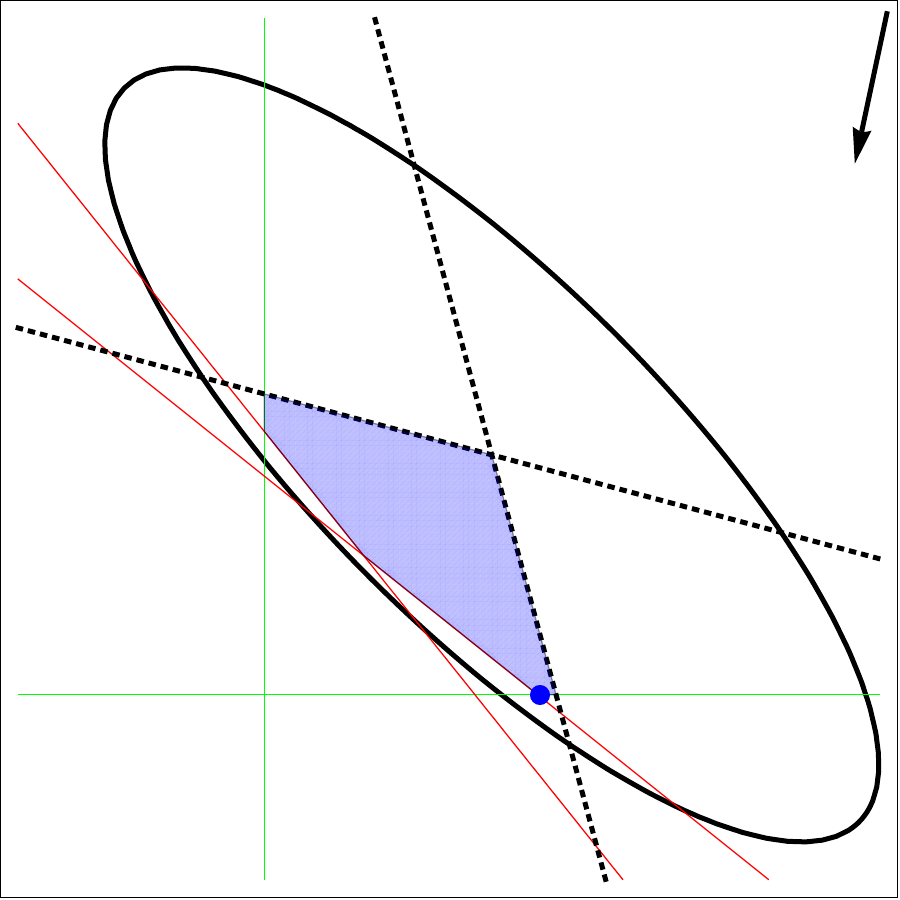}
	}
	\subfigure[SOCP is inexact due to reactive and active lower bounds.]{
            \includegraphics[width=0.465\columnwidth]{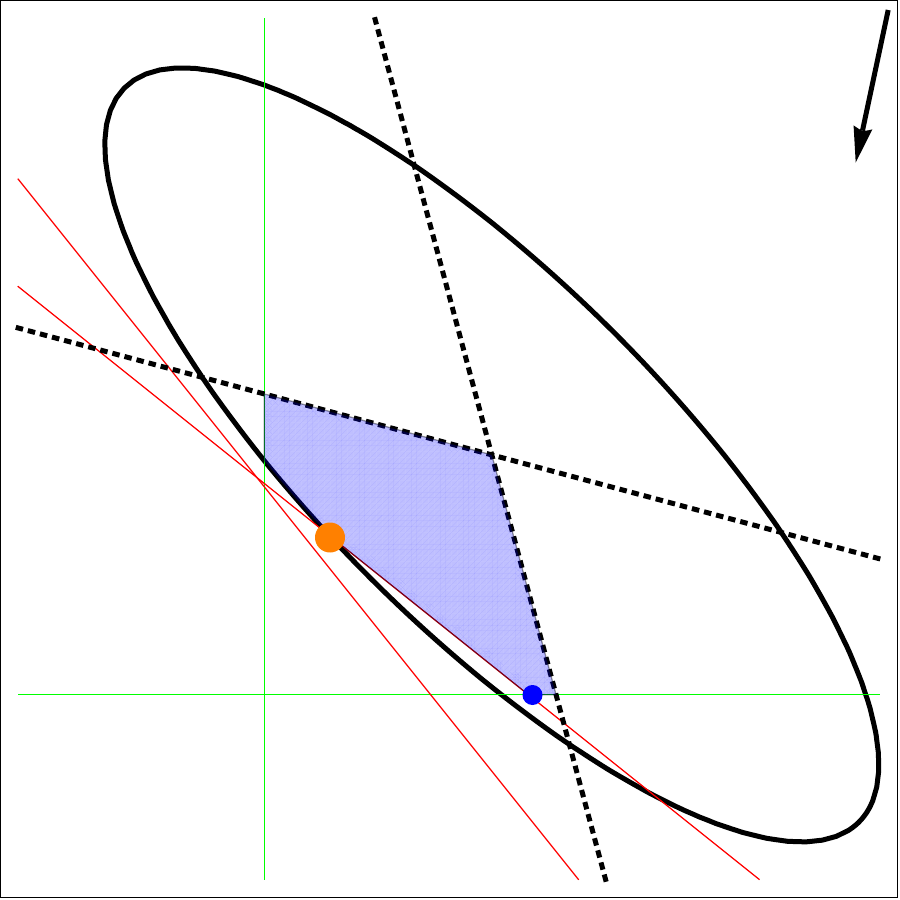}
	}
	\subfigure[SOCP is inexact due to angle and {reactive} lower bounds.]{
            \includegraphics[width=0.465\columnwidth]{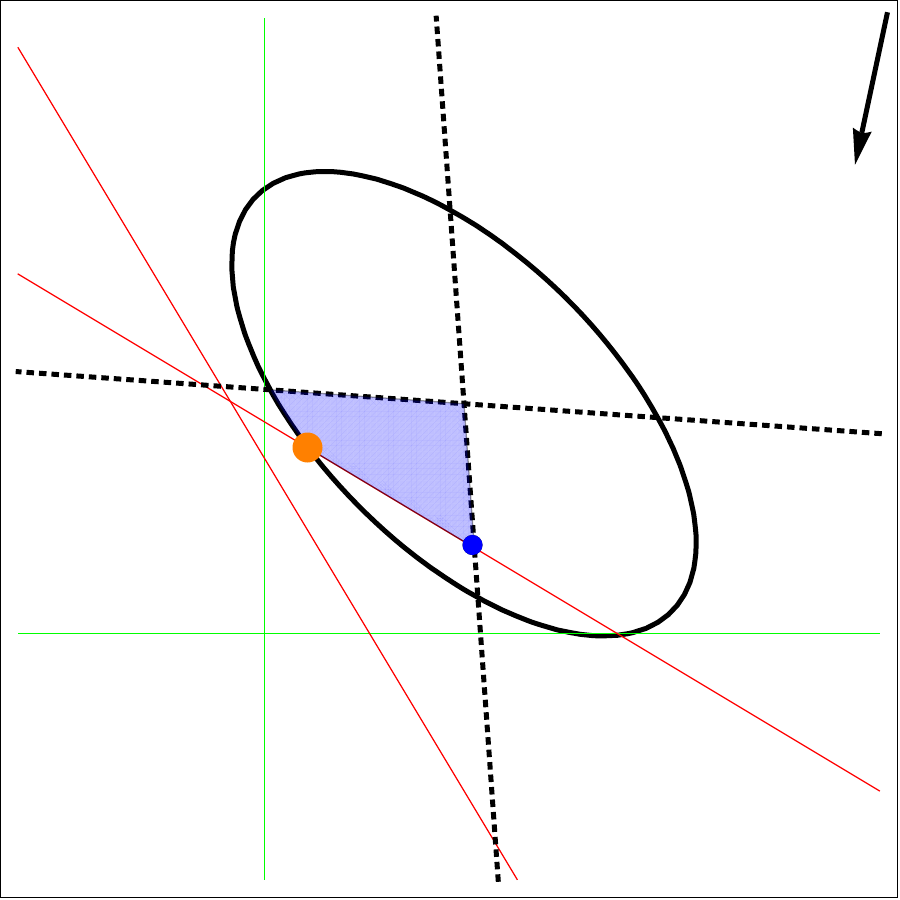}
         }
\end{center}
\end{figure*}

If $c_{22}^E <  \underline c_{22}$, where $\underline c_{ii}:=(V_i^{\min})^2$, then SOCP is feasible while OPF is infeasible. An example of this case can be seen from Figure \ref{Cases}b, which shows the zoomed in part of the hyperbola.

\item If $c_{11}^O - c_{22}^O < \Delta$, $c_{22}^E \ge  \underline c_{22}$, and $c_{11}^E\geq \underline c_{11}$,
then the SOCP relaxation is exact as in Figure \ref{Cases}c. In fact, any point in the convex combination of $c^R$ and $c^E$ is SOCP optimal. Such a point can always be {corrected} by reducing $c_{11}$, $c_{22}$ and $c_{12}$ components by the same amount until we reach $c^E$, which is the OPF optimal solution.

\item If $c_{11}^O - c_{22}^O < \Delta$, $c_{22}^E \ge  \underline c_{22}$, and $c_{11}^E  < \underline c_{11}$, define the intersection  of $c_{11} - c_{22} = \Delta$ with the bounding lower bound of either $c_{11}$ or $c_{22}$:
\begin{equation}
\hspace{0mm}(c_{11}^L, c_{22}^L) = \begin{cases}   ( \underline c_{22} +  \Delta, \underline c_{22} ) \text{ if }   \underline c_{11} -\underline c_{22} \le \Delta \\
(\underline c_{11}, \underline c_{11} -\Delta) \text{ o.w. }
\end{cases}
\end{equation}
Observe that any point in the convex combination of $c^L$ and $c^R$ is SOCP optimal.  However, there is no feasible OPF solution with the same objective function value.
Lastly, let us define the lower intersection  of the hyperbola (\ref{feasCurve}) and the $c_{11}$ lower bound as $(c_{11}^I, c_{22}^I)$, where $c_{11}^I=\underline c_{11}$, and $c_{22}^I$ as
\begin{equation}
c_{22}^I = \frac{2\beta + \underline c_{11} - \sqrt{(2\beta + \underline c_{11})^2 - 4(\alpha^2+\beta^2)}}{2}\label{eq:cI}
\end{equation}
We have two final cases:
\begin{itemize}
\item If $c_{22}^I <  \underline c_{22}$, any point in the convex combination of $c^L$ and $c^R$ is SOCP optimal.  However,  OPF is infeasible. An example of this case can be seen in  Figure \ref{Cases}d.
\item If $c_{22}^I \ge  \underline c_{22}$, any point in the convex combination of $c^L$ and $c^R$ is SOCP optimal.  However, OPF has a unique optimal solution at $c^I $ as can be seen  in  Figure \ref{Cases}e. Hence, relaxation is inexact. Assuming a linear cost function with  coefficient $1$, optimality gap can be calculated as $-G(c_{22}^L-c_{22}^I)$.
\end{itemize}
\end{itemize}

The above analysis proves the following theorem.
\begin{thm}
In a two-bus one-generator system {with linear objective}, the SOCP/SDP relaxation of the AC OPF problem has the following possible outcomes:
\begin{enumerate}
\item SOCP relaxation is exact: If  $c_{11}^O - c_{22}^O \ge \Delta$ or  if $c_{11}^O - c_{22}^O < \Delta$, $c_{22}^E \ge  \underline c_{22}$, $c_{11}^E \ge \underline c_{11}$.
\item SOCP relaxation is inexact with finite optimality gap: If  $c_{11}^O - c_{22}^O < \Delta$, $c_{22}^E \ge  \underline c_{22}$,  $c_{11}^E  < \underline c_{11}$, $c_{22}^I \ge  \underline c_{22}$. The optimality gap is $-G(c_{22}^L-c_{22}^I)$.
\item {SOCP relaxation is feasible and OPF is infeasible}: If $c_{11}^O - c_{22}^O < \Delta$, $c_{22}^E <  \underline c_{22}$ or if $c_{11}^O - c_{22}^O < \Delta$, $c_{22}^E \ge  \underline c_{22}$,  $c_{11}^E  < \underline c_{11}$, $c_{22}^I <  \underline c_{22}$.
\end{enumerate}
Here, $c^O, c^E, c^R, c^L, c^I$ are defined in \eqref{eq:cO}-\eqref{eq:cI}, respectively.
\end{thm}
%
%
%
%
%
%


\section{Examples of Inexact SOCP Relaxations} \label{section:examples}
We have obtained a complete characterization for a 2-bus network with a single generator, and shown that the SOCP relaxation is exact only under certain conditions. In this section,
we present further counterexamples of radial networks with two and three buses. Most of the network parameters are selected from IEEE test instances. Transmission line capacity is assumed to be large. For all the buses, $V_i^\text{min}=0.9$ and  $V_i^\text{max}=1.1$. Production costs are taken as linear functions. OPF problem with alternative  formulation (\ref{SOCP}) is solved to global optimality with BARON \cite{BARON}. SOCP relaxations are solved using interior point solver MOSEK \cite{Mosek}.
\vspace{-3mm}
\subsection{2-Bus, 2-Generator Example}
Let us consider a 2-bus network with one generator located at each bus. Data of this example is  given in Table \ref{BusGen2}. The impedance of line (1,2) is   $0.01008+\mathrm{i} 0.0504$.

\begin{table}[h!]
\caption{Bus and generator data for 2-bus 2-generator example.}
\label{BusGen2}
\begin{center}
\begin{tabular}{rrrrrrrr}
\hline
Bus &         $p_i^d$ &         $q_i^d$  &     $p_i^{\text{min}}$ & $p_i^{\text{max}}$     &      $q_i^{\text{min}}$ & $q_i^{\text{max}}$  &       cost \\
\hline
         $1$ &         $75$ &      $-84.7$ &         $75$ &        $250$ &        $-30$ &        $300$ &          $5.0$ \\
         $2$ &        $105$ &       $22.8$ &         $70$ &        $300$ &        $-30$ &        $300$ &        $1.2$ \\
\hline
\end{tabular}
\end{center}
\end{table}


\begin{figure*}
\begin{center}
\caption{Projection of feasible region of 3-bus example onto $(p_1^g,q_1^g)$ space with respect to different load levels. Horizontal axis is $p_1^g$ and vertical axis is $q_1^g$. Black curve and blue region are the feasible regions of OPF and SOCP relaxation, respectively.  Red line is the lower bound on $q_1^g$. All figures are in p.u.}  \label{FeasReg3}
	\subfigure[$\gamma=0.90$]{
            \includegraphics[width=0.4645\columnwidth]{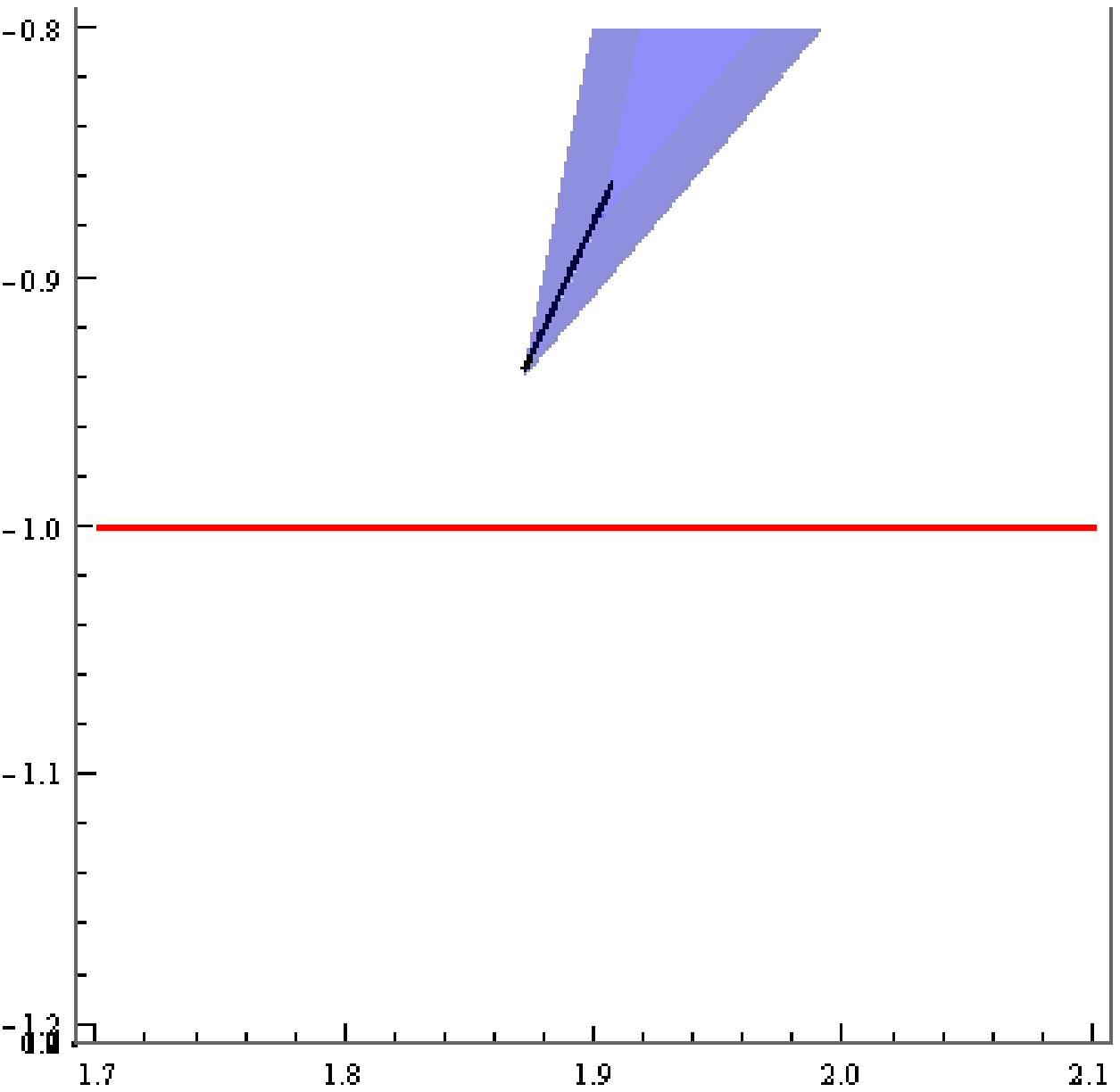}
	}
	\subfigure[$\gamma=1.00$]{
            \includegraphics[width=0.4645\columnwidth]{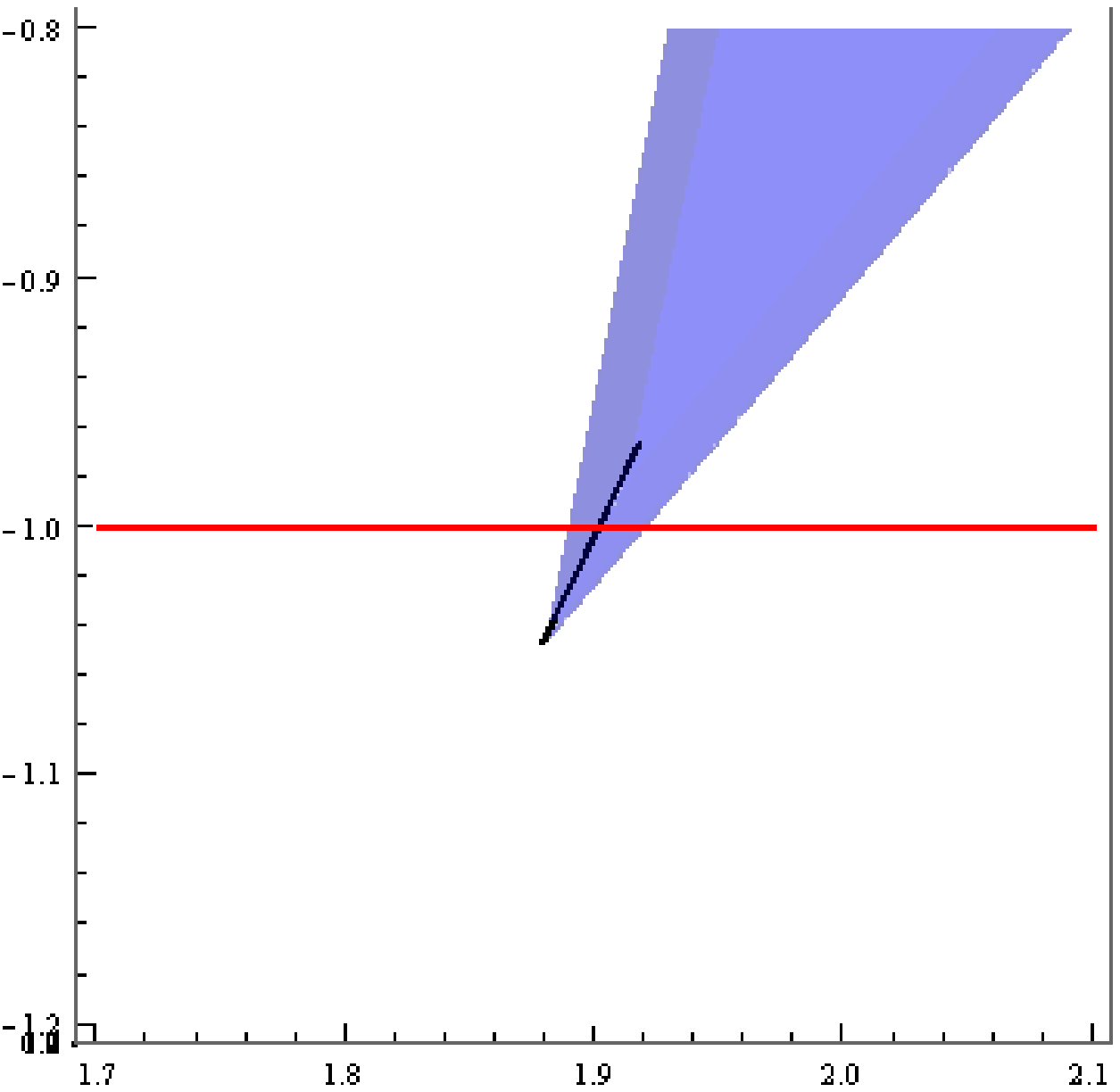}
	}
	\subfigure[$\gamma=1.10$]{
            \includegraphics[width=0.4645\columnwidth]{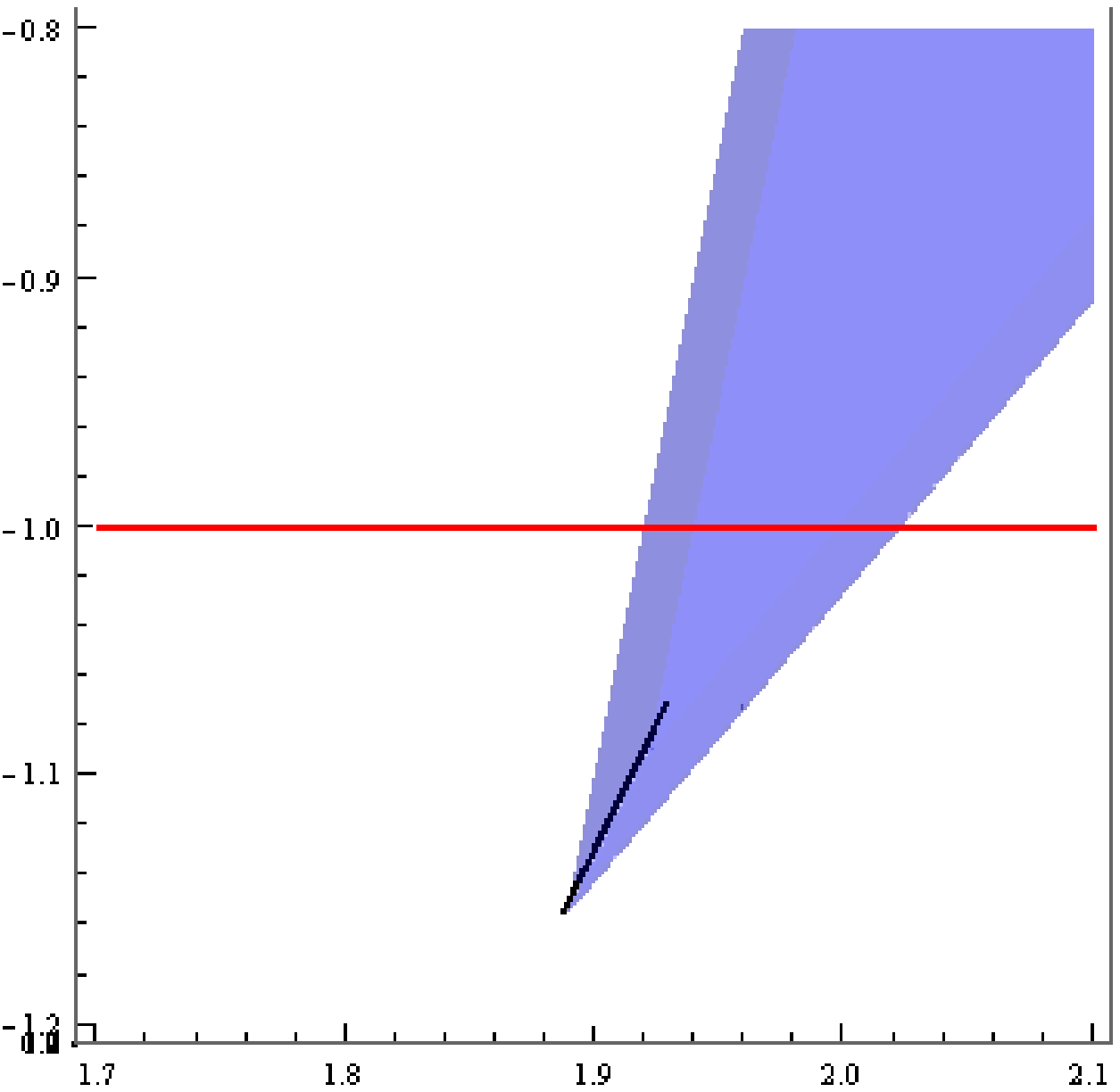}
           }
\end{center}
\end{figure*}

In Table \ref{table:Sens2}, we compare the SOCP relaxation and the global optimal solution of OPF for different levels of load, where load is varied as $[p_1^d\;\; p_2^d \;\; q_1^d \;\; q_2^d] = \gamma\cdot[ 75 \;\; 105 \;\; -84.7  \;\; 22.8 ]$ for some positive parameter $\gamma$.

\begin{table}[h!]
\caption{Objective costs for 2-bus 2-generator with varying load.}
\label{table:Sens2}
\begin{center}
\begin{tabular}{rrr}
\hline
     $\gamma$ &        OPF &       SOCP  \\
\hline

     $0.12$ &   infeasible   &   infeasible    \\

     $0.13$ &   infeasible   &    $459.00$     \\

     $0.80$ &   infeasible   &   $459.00$     \\

     $0.81$ &    $460.13$     &    $460.13$   \\

      $0.98$ &   $496.96$   &   $496.96$   \\

      $0.99$ &  $499.15$   &  $499.15$   \\

      $1.00$ & $563.56$  &   $501.46$ \\
     $1.01$ &   $641.21$  &     $503.76$  \\
     $1.02$ &   infeasible &   $506.07$   \\
     $2.92$ &   infeasible &      $1608.75$    \\
     $2.93$ &   infeasible &  infeasible   \\
\hline
\end{tabular}
\end{center}
\end{table}

When $\gamma \in [0.81, 0.99]$, we observe that the SOCP relaxation is exact. For $\gamma$ around 1.00, there is a finite optimality gap, which can be as large as $21.44\%$ at $\gamma=1.01$. Finally, for $\gamma \ge 1.02$, OPF becomes infeasible, whereas SOCP relaxation is still feasible. In fact, SOCP relaxation fails to detect infeasibility of OPF problem until $\gamma$ exceeds $2.93$.

Now, let us consider the case where voltages are fixed. In  \cite{Lavaei14}, it has been proven that if angle differences are guaranteed to be small enough, then SDP/SOCP relaxations are tight even if there are real power lower bounds. However, we present an example which demonstrates that this does not extend to the case with reactive power lower bounds. To this end, let us fix the squared voltage magnitudes to  $(c_{11}, c_{22})=(0.874,0.816)$. In this case, the global optimal solution of OPF is $573.82$ while the SOCP relaxation gives $503.37$. Hence, there is an optimality gap, even though angle difference is less than  $1^{\circ}$.

%

Figure \ref{FeasReg2 all} presents possible configurations of the feasible region of the OPF and the SOCP relaxation projected to the $(p_1^g, p_2^g)$ space. In Figure \ref{FeasReg2 all}a, the SOCP relaxation is exact, while  in Figure \ref{FeasReg2 all}b, the OPF is infeasible although  the SOCP is feasible. In Figure \ref{FeasReg2 all}c, the SOCP relaxation is inexact due to the combined effect of active and reactive lower bounds. Finally, in Figure \ref{FeasReg2 all}d, the SOCP relaxation is inexact due to practical angle bounds ($30^\circ$) and reactive lower bounds.

\subsection{3-Bus, 1-Generator Example}

Consider a 3-bus radial network with three loads $[p_1^d\;\; p_2^d \;\; p_3^d] = [ 50 \;\; 70 \;\; 60]$ and $[q_1^d \;\; q_2^d \;\; q_3^d] = [ -52.3 \;\; 14.1 \;\;  -82.3]$. The impedance of lines (1,2)  and (2,3) are    $0.01008+\mathrm{i} 0.0504$ and $0.07500+\mathrm{i} 0.0840$.
The only generator is located at bus 1 with $150  \le p_1^g \le 550$ and $ -100  \le q_1^g \le 500$. The cost of power generation is \$5 per MW.
Assume that the reactive load is scaled as $[q_1^d \;\; q_2^d \;\; q_3^d] = \gamma[ -52.3 \;\; 14.1 \;\;  -82.3]$ for some positive $\gamma$. Table \ref{table:Sens3} shows the optimal costs of the OPF and the SOCP relaxation for different values of $\gamma$.

\begin{table}[h!]
\caption{Objective costs for 3-bus example with varying load.}
\label{table:Sens3}
\begin{center}
\begin{tabular}{rrr}
\hline
     $\gamma$ &        OPF &       SOCP  \\
\hline
      $0.95$ & $939.45$ &   $939.45$ \\

      $0.96$ & $939.90$ & $939.90$ \\

      $0.97$	& $941.57$	& $940.87$ \\

      $1.00$ & $950.70$ &   $945.45$ \\
     $1.03$ & $959.91$ &   $950.05$\\
     $1.04$ &   infeasible &  $951.60$ \\
\hline
\end{tabular}
\end{center}
\end{table}

For small values of $\gamma$, e.g. $\gamma\leq 0.96$, SOCP is exact. For values around $\gamma=1$, we observe a finite optimality gap between OPF and SOCP, where for larger values of $\gamma\geq 1.04$, the OPF becomes infeasible while the SOCP relaxation is still feasible. The infeasibility is exactly caused by the lower bound on reactive generation power.

For this example, we also give the feasible region of OPF problem projected onto the $(p_1^g,q_1^g)$ space in Figure \ref{FeasReg3}. When $\gamma=0.90$, reactive power lower bound is redundant and the optimal solution of SOCP relaxation is feasible for OPF. However, for  $\gamma=1.00$, constraint $q_1^{\text{min}}  \le q_1^g$ is binding. Note that the optimal solution of the SOCP relaxation is not feasible for OPF and hence, the relaxation is not exact.  Finally, when $\gamma=1.10$, SOCP is feasible whereas OPF is infeasible.

\section{Library of Radial Networks with Inexact SDP/SOCP Relaxation}
\label{section:library}

\subsection{Generation of Instances}
To facilitate further research, we generate several radial network instances from meshed networks in  MATPOWER  \cite{Matpower}. Given a meshed network, we first find a spanning tree { by switching off lines to obtain a radial network}. Then, only load values and generation lower/upper bounds are changed, whenever necessary, to construct examples where the SOCP/SDP relaxation is not exact. New instances can be downloaded from {\url{https://sites.google.com/site/burakkocuk/}}.

{Our examples are based on 9-, 14-, 30-, 39- and 57-bus standard instances. Due to our construction of the network topology, AC feasibility becomes a major issue. Although unrealistic examples can be constructed for even larger networks by reducing load values considerably, we choose not to sacrifice the realistic features of the instances. }


\subsection{Computational Results for  SDP Relaxation vs. Global Optimal Solution}
For each instance generated as described above, we solve the SDP relaxation using MOSEK  \cite{Mosek}. The code is written in C\# language and Visual Studio 2010 is used as the compiler.  We report the value of the objective function, computation time and the rank of the solution. Here, rank is determined as the number of eigenvalues that are larger than $10^{-5}$.

SDP relaxation is compared   against global optimal solution found using BARON \cite{BARON} and local solution found by MATPOWER \cite{Matpower} { and IPOPT \cite{wachter}}. Relative optimality gap for BARON is set to 0 so that global optimality can be certified. We should note that performance of BARON on rectangular formulation (\ref{rect form}) is very poor as it requires hours to prove global optimality. Instead, we use reformulation (\ref{SOCP}), which is valid for radial networks.

\begin{landscape}
\begin{table*}
\caption{SDP relaxation  vs. global solver BARON and local solver MATPOWER.}
\label{SDP vs. BARON}
\begin{center}
\begin{tabular}{ccrrrrrrrrrrrr}
\hline
           &            &  \multicolumn{ 3}{c}{SDP Relaxation} & & \multicolumn{ 2}{c}{BARON} &                      & \multicolumn{ 2}{c}{MATPOWER} &           & \multicolumn{ 2}{c}{  {IPOPT} } \\
\cline{3-5} \cline{7-8} \cline{10-11}\cline{13-14}
  based on &       type &  objective &       time(s) &       rank & & objective &       time(s) &             \% gap    &objective &       time(s)  & &objective &  time(s)\\
\hline
     case9 &  quadratic &    5335.70 &       0.04 &          8 & &  11277.95 &       1.17 &         52.69          &      - &       0.17 &  & -	&0.17 \\

     case9 &     linear &    1481.93 &       0.06 &          8 &  &  1756.47 &       1.11 &           15.63                &  - &       0.08&&   -	&0.20\\

    case9Q &  quadratic &   10835.70 &       0.04 &          8 & &  16778.87 &       1.36 &          35.42           &      - &       0.08&&16779.48&	0.31 \\

    case14 &  quadratic &   11861.87 &       0.07 &          8 & &  11932.07 &      35.32 &             0.59             & 11932.25 &       0.11 &&11932.25&	0.28\\

    case14 &     linear &    9892.70 &       0.09 &          4 &  &  9952.42 &       0.79 &             0.60                 & 9952.59 &       0.09&&9952.58&	0.23 \\

case\_ieee30 &  quadratic &    4244.53 &       0.17 &         12 & &   4336.03 &    8347.79 &      2.11  	 &        - &       0.12&&4794.32&	0.15 \\

case\_ieee30 &     linear &    3035.61 &       0.22 &         12 &  &  3606.91 &    2494.31 &     15.84     	 &     - &       0.09 &&4562.26	&0.14\\

    case30 &  quadratic &     607.72 &       0.15 &          8 & &    619.01 &       2.52 &            1.82    		& 619.04 &       0.09&&619.04	&0.23 \\

    case30 &     linear &     435.58 &       0.23 &          6 &  &   445.83 &       8.50 &           2.30    		& 445.84 &       0.11 &&445.84	&0.14\\

   case30Q &  quadratic &     676.88 &       0.20 &          4 &  &   690.06 &       5.16 &             1.91  		 &  690.08 &       0.11 &&690.08&	0.36\\

    case39 &  quadratic &   44869.01 &       0.29 &          4 & &  45035.32 &     110.59 &        0.37 	 &         - &       0.14 &&45037.05&	0.27\\

    case39 &     linear &    1900.09 &       0.36 &          4 &  &  1903.07 &    1566.88 &            0.16       	 &  - &       0.16&& 1903.14	&0.15\\

  {  case57 }  &  quadratic &   10458.06 &       0.92 &          20 & & 12100.00 &  $>10800$ &    13.57   	  &  12100.90 &       0.15 &&12100.86	&0.27\\

  {  case57 } &     linear &    8399.82 &       0.96 &          20 &  &  10173.10 & $>10800$  &   17.43 		 &    10173.00&       0.16 &&10172.98&	0.26\\
\hline
\end{tabular}

\end{center}
\end{table*}

\end{landscape}

For all experiments, we used a 64-bit computer with Intel Core i5 CPU 3.33GHz processor and 4 GB RAM. Each instance is solved twice with quadratic and linear objectives. For the latter, we simply ignore the quadratic cost coefficients.

Our findings are summarized in  Table \ref{SDP vs. BARON}. One can  see that the SDP relaxation solution can be of high-rank  (up to 12 for case\_ieee30 {and 20 for case57}). Also, the optimality gap (column ``\% gap'') {computed as $100\times (1 - z^{  SDP} / z^{ BARON})$, where $z^{ SDP}$ and $ z^{ BARON}$ are respectively the values of the SDP relaxation and the global optimal solution found by  BARON,} can be quite large (more than 52\% for case9 with quadratic objective). Our examples clearly show that {the optimal value of the SDP relaxation can be quite different from the global optimal value}. {We also compare the optimal dispatch solutions $p^{SDP}$ and $p^{BARON}$ computed by the SDP relaxation and BARON to show that large differences in the objective function values are not  artifacts of the cost parameters. In fact, the 2-norm $\|p^{SDP} - p^{BARON}\|$ is  large, varying from 0.16 p.u. to 3.16 p.u. for our instances. This illustrates that the optimal solutions are quite different from one another.}

In general, MATPOWER is accepted to be a reliable and efficient OPF solver.  It manages to  find the global optimal solution up to a negligible difference for seven of the instances from our library. However, we observe that it fails to solve the remaining seven instances due to numerical issues. 
{There are other robust NLP solvers available, e.g. IPOPT, which gives near global optimal solution for nine instances in the library, where small discrepancies in optimal objective function values compared to BARON are due to numerical errors. On the other hand, IPOPT fails in two 9-bus examples and it finds suboptimal solutions three times for both of the case\_ieee30 instances and case39 with quadratic objective.}

We should note that the global solver BARON can be computationally expensive. For instance, for case\_ieee30 with a quadratic objective, it requires more than 2 hours to prove optimality {whereas for 57-bus instances, BARON  is not able to certify the global optimal solution within 3 hours time limit. Upon termination, the optimality gaps are  38.40\%  and 29.37\% for quadratic and linear objectives, respectively.}
Also, the reformulation of OPF \eqref{SOCP} is only valid for radial networks. Hence, in general, using BARON as it is may not be applicable to
large-scale OPFs.


\section{Bound Tightening and Valid Inequalities for Global Optimization}
\label{section:valid}

In this section, we propose valid inequalities for the SOCP relaxation of the OPF problem to improve the computational time of the global solver BARON. The main algorithm of BARON is based on  spatial branch-and-bound \cite{BARON}. It utilizes convex envelopes of the feasible region and polyhedral relaxations to improve lower bounds and prove global optimality. Therefore, it is very important to add valid inequalities and variable bounds so that BARON can obtain tighter relaxations. {A more detailed description is provided in Appendix \ref{app:baron}.}

To begin with, let us focus on formulation (\ref{SOCP}). Observe that $c_{ij}$ and $s_{ij}$  do not have explicit variable bounds although they have implied bounds due to (\ref{voltageAtBusR})  and (\ref{coupling}) as
\begin{equation}
-V_i^{\text{max}}V_j^{\text{max}} \le c_{ij}, s_{ij} \le V_i^{\text{max}}V_j^{\text{max}} \quad (i,j) \in \mathcal{L}
\end{equation}
However, these bounds are very loose knowing that angle differences are generally small. This fact  suggests that these bounds can be improved. One way to obtain variable bounds is to optimize $c_{ij}$ and $s_{ij}$  over the set $\mathcal{S} = \{(p,q,c,s) : (\ref{activeAtBusR}) - (\ref{coupling}) \}$, which is a nonconvex set. Alternatively, one can find weaker bounds over the set $\mathcal{S}' = \{(p,q,c,s) : (\ref{activeAtBusR}) - (\ref{sine}),  (\ref{couplingSOCP})  \}$ by solving SOCP relaxations. Let  $\underline c_{ij} $ ($\underline s_{ij}$) and $\overline c_{ij} $ ($\overline s_{ij} $)  denote  lower and upper bounds found for $c_{ij}$ ($s_{ij}$), respectively.

Now, let us investigate how the box $\mathcal{B}_{ij}=[\underline c_{ij}, \overline c_{ij}]\times[\underline s_{ij}, \overline s_{ij}]$  is positioned with respect to the {``ring"-like set} $\mathcal{R}_{ij} = \{(c_{ij},s_{ij}): {\underline R_{ij}^2} \le c_{ij}^2 + s_{ij}^2 \le {\overline R_{ij}^2} \}$ where $\underline R_{ij}=V_i^{\text{min}}V_j^{\text{min}}$ and $\overline R_{ij}=V_i^{\text{max}}V_j^{\text{max}}$. In our experiments, {we observe that $\underline c_{ij} > 0$, which we assume hereafter. We should note that this is not a restrictive assumption, similar valid inequalities described below can be generated even if this assumption does  not hold.}

\begin{figure}
\begin{center}
\caption{Positioning of  $\mathcal{B}_{ij}$ and  $\mathcal{R}_{ij}$. Red line is the cut produced by Algorithm \ref{alg:T1C}, when applicable. }  \label{typical}
	\subfigure[Case 1]{
            \includegraphics[width=0.373\columnwidth]{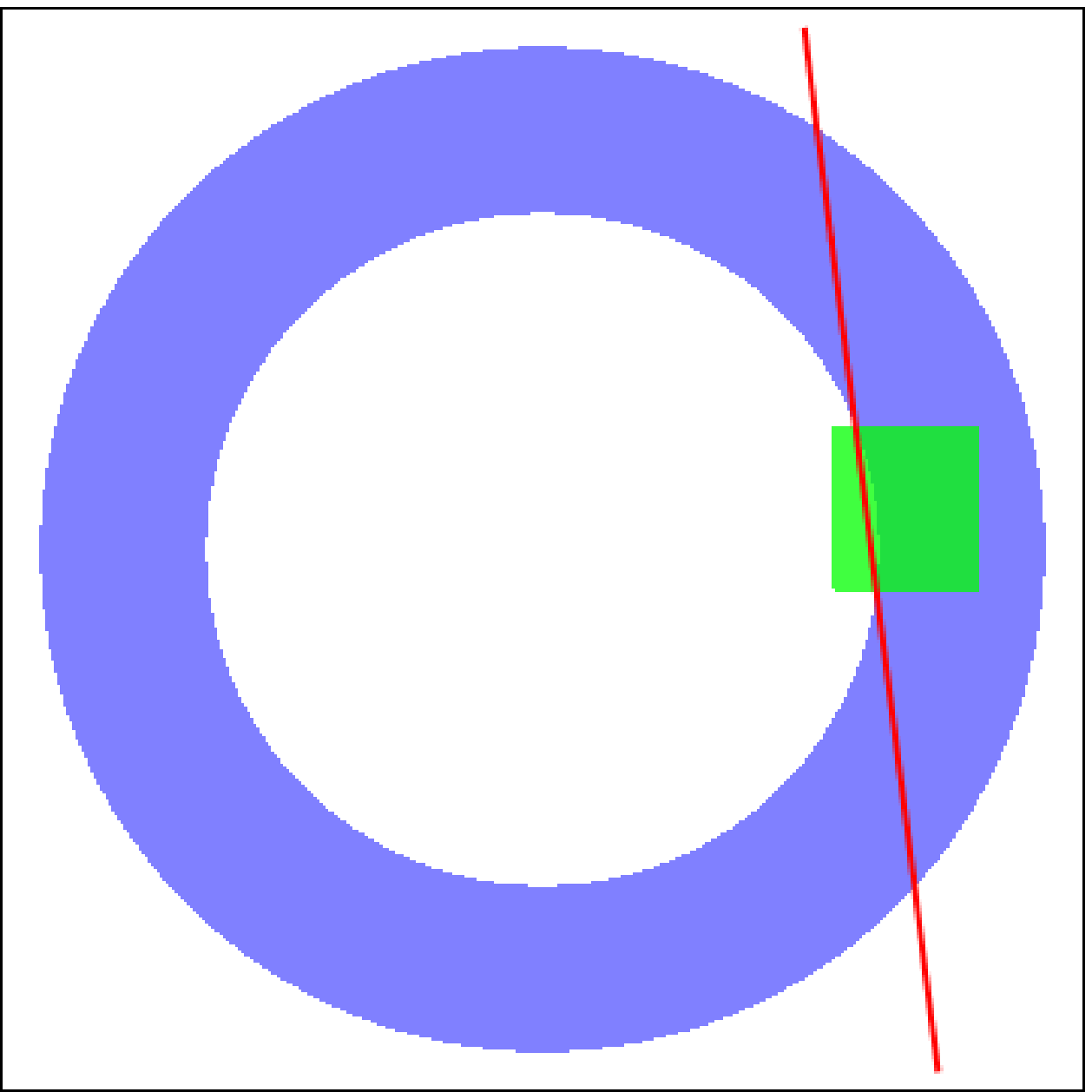}
	}
	\subfigure[Case 2]{
            \includegraphics[width=0.373\columnwidth]{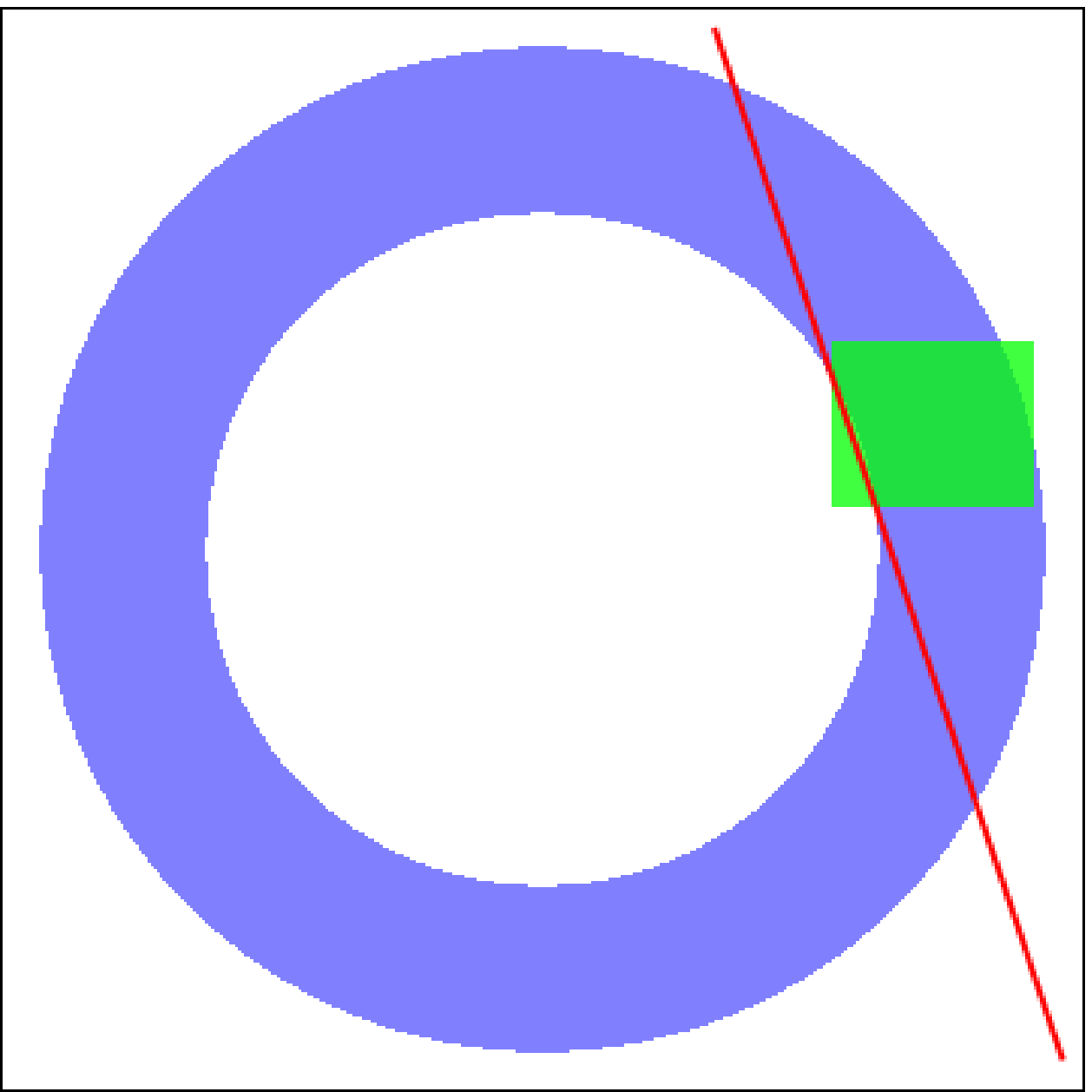}
	}
	\subfigure[Case 3]{
            \includegraphics[width=0.373\columnwidth]{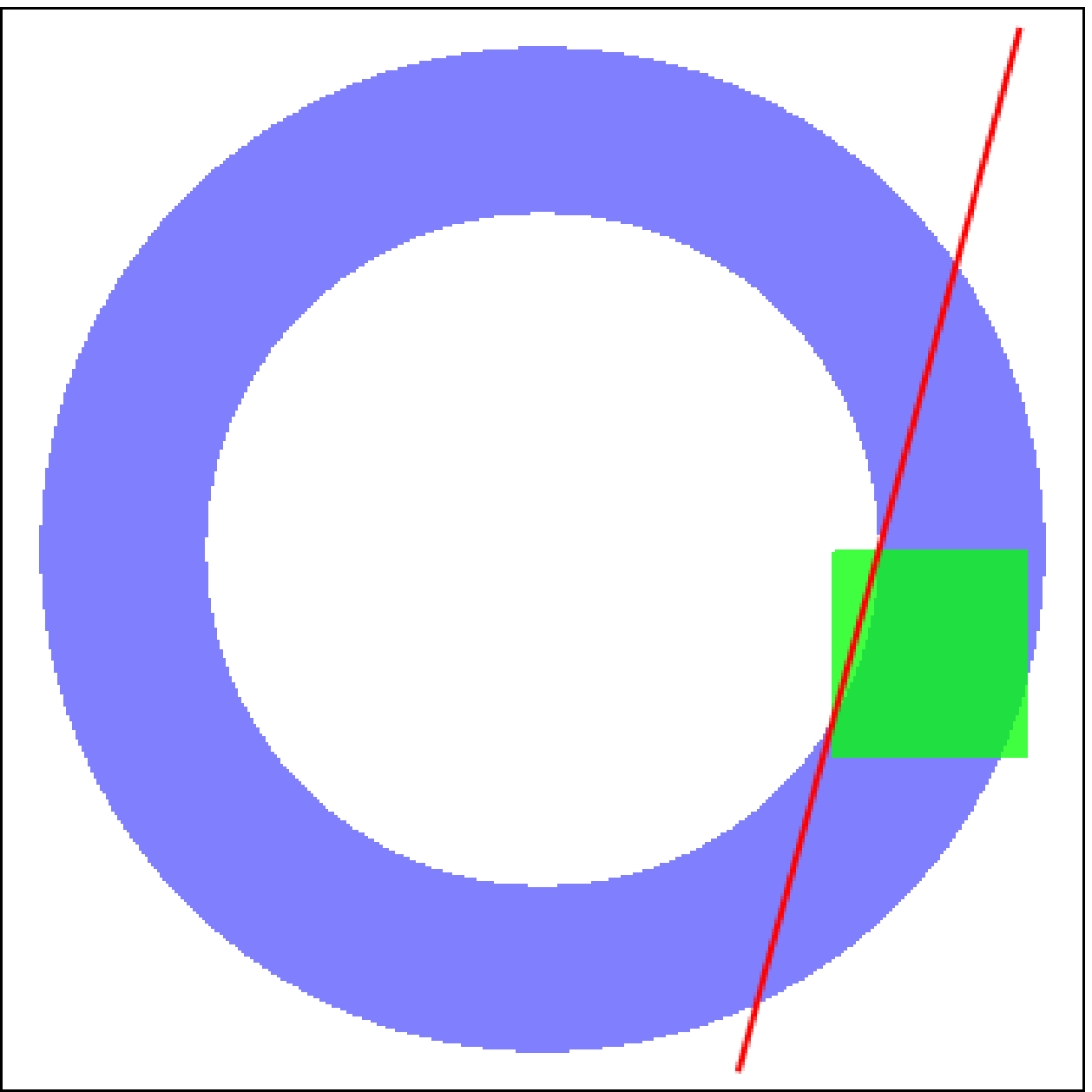}
	}
	\subfigure[Case 4]{
            \includegraphics[width=0.373\columnwidth]{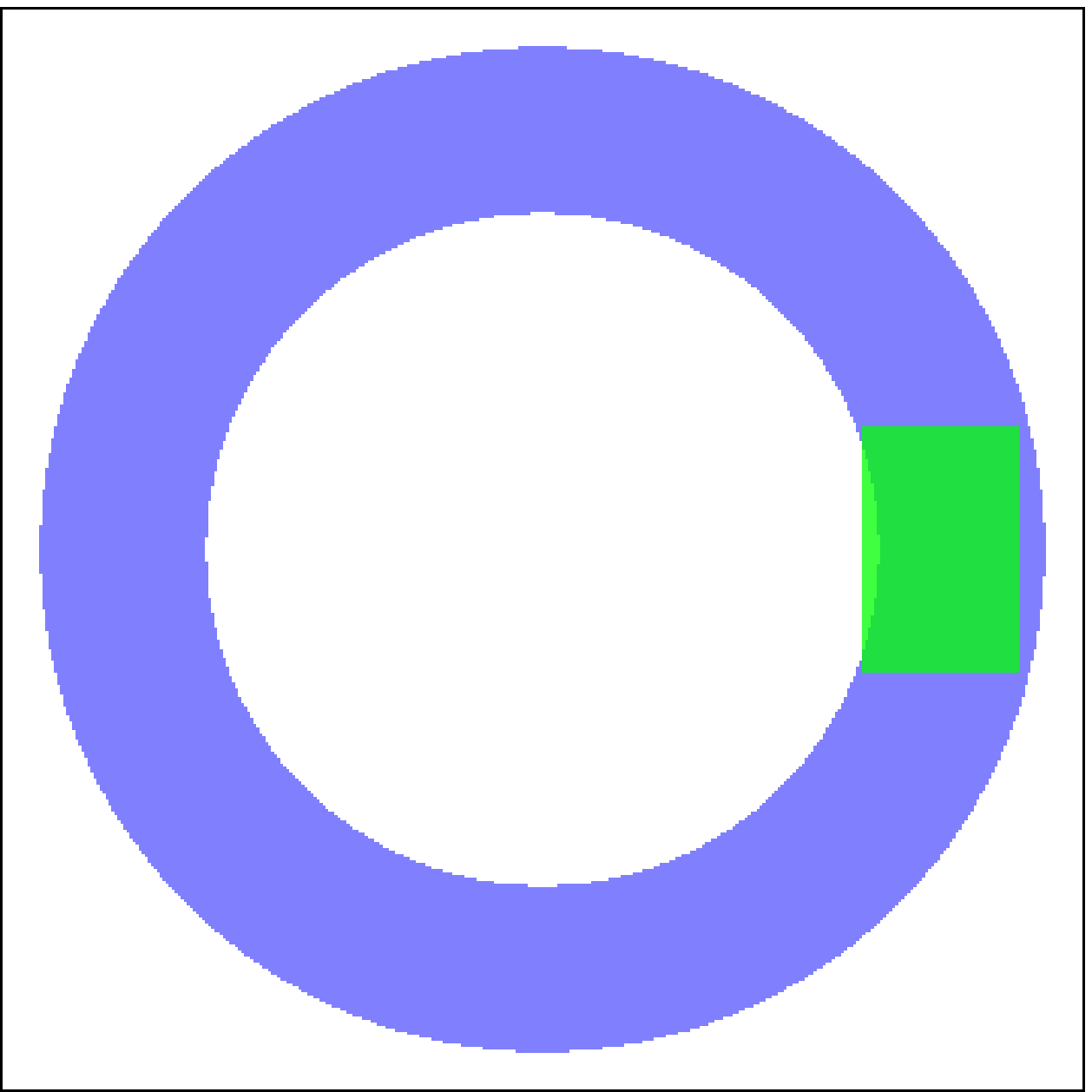}
	}
\end{center}
\end{figure}

 Let us focus on the case with  $\underline c_{ij} < \underline R_{ij}$, which gives rise to four possibilities:
\begin{itemize}
\item Case 1: $ \|(\underline c_{ij}, \underline s_{ij})\| < \underline R_{ij}$, $ \|(\underline c_{ij}, \overline s_{ij})\|< \underline R_{ij}$
\item Case 2: $ \| (\underline c_{ij}, \underline s_{ij})\| < \underline R_{ij}$, $ \|(\underline c_{ij}, \overline s_{ij})\|\ge \underline R_{ij}$
\item Case 3: $ \| (\underline c_{ij}, \underline s_{ij})\| \ge \underline R_{ij}$, $\|(\underline c_{ij}, \overline s_{ij})\|< \underline R_{ij}$
\item Case 4: $ \| (\underline c_{ij}, \underline s_{ij})\| \ge \underline R_{ij}$, $\|(\underline c_{ij}, \overline s_{ij})\|\ge \underline R_{ij}$
\end{itemize}
Figure \ref{typical} shows typical examples for each of four cases.  In the rest of this section, we concentrate on how we can obtain valid inequalities for Cases 1, 2, and 3.

\subsection{Valid Inequalities}
These cuts are designed to cut off the portion of $\mathcal{B}_{ij}$ inside the inner circle for Cases 1, 2, and 3 as depicted in Figure \ref{typical}.  Algorithm \ref{alg:T1C} gives the exact procedure. {Note that the validity of the inequality follows from the fact that the points cut off from the box  have norm less than $\underline R_{ij}$.} Note that for Case 4, the algorithm would produce the inequality $c_{ij} \ge \underline c_{ij}$, hence it is omitted.


\begin{algorithm}
\caption{Generation of Valid Inequalities.}
\label{alg:T1C}
\begin{algorithmic}
 \FORALL{$(i,j) \in \mathcal{L}$}
\STATE Compute $\underline c_{ij}$,  $\overline c_{ij}$,  $\underline s_{ij}$ and $\overline s_{ij}$ over $ \mathcal{S}'$.
\STATE Update $ \mathcal{S}' = \mathcal{S}' \cap \mathcal{B}_{ij}$.

\IF{ ${0  < } \underline c_{ij} < \underline R_{ij}$ }
	\IF{       $\|(\underline c_{ij}, \underline s_{ij})\| < \underline R_{ij}$, $ \|(\underline c_{ij}, \overline s_{ij})\|< \underline R_{ij}$ }
		\STATE  Set $y_1=\overline s_{ij}$, $y_2=\underline s_{ij}$ and compute \\ $x_1=\sqrt{\underline R_{ij}^2-\overline s_{ij}^2}$, $x_2=\sqrt{\underline R_{ij}^2-\underline s_{ij}^2}$
	\ELSIF{  $\|(\underline c_{ij}, \underline s_{ij})\| < \underline R_{ij}$, $\|(\underline c_{ij}, \overline s_{ij})\|\ge \underline R_{ij}$}
		\STATE   Set $x_1=\underline c_{ij}$, $y_2=\underline s_{ij}$ and compute \\ $y_1=\sqrt{\underline R_{ij}^2-\underline c_{ij}^2}$, $x_2=\sqrt{\underline R_{ij}^2-\underline s_{ij}^2}$
	\ELSIF{  $\| (\underline c_{ij}, \underline s_{ij})\| \ge \underline R_{ij}$, $\|(\underline c_{ij}, \overline s_{ij})\|< \underline R_{ij}$ }
		\STATE Set $y_1=\overline s_{ij}$, $x_2=\underline c_{ij}$  and compute \\ $x_1=\sqrt{\underline R_{ij}^2-\overline s_{ij}^2}$, $y_2=-\sqrt{\underline R_{ij}^2-\underline c_{ij}^2}$
	\ENDIF
\STATE  Add  $(y_1-y_2)c_{ij}-(x_1-x_2)s_{ij} \ge x_2y_1-x_1y_2$ as a valid inequality and update $ \mathcal{S}'$.
\ENDIF
 \ENDFOR
\end{algorithmic}
\end{algorithm}

\subsection{Numerical Experiments}

The effect of valid inequalities are tested on our library of instances. The results are summarized in Table \ref{BARON cut}. We should note that  MATPOWER is very efficient and accurate for the seven instances it is able to solve as shown in Table \ref{SDP vs. BARON}. Therefore, we mainly focus on the other seven instances where MATPOWER fails to solve. In Table \ref{BARON cut}, preprocessing refers to computing variable bounds and valid cuts.

For the 9-bus instances (case9, case9Q), BARON's computation time reduces slightly with the addition of cuts. However, the preprocessing time dominates the total computation time, which is larger than the case without cuts.

For the 30-bus IEEE instances, BARON can require hours to terminate. With the addition of variable bounds, total computation time reduces by 89\% and 90\% for quadratic and linear objectives, respectively. Quite impressively, the inclusion of valid inequalities further reduces the total computation time to only 17 seconds, less than 0.1\% of the computation time without variable bounds and cuts.

For 39-bus instances, the addition of variable bounds brings down total computation time by 76\% and 95\%  for quadratic and linear objectives, respectively. In this case, the inclusion of valid inequalities decreases the computational time  for linear objective. On the other hand, cuts slightly increases the total computational time  in the case of quadratic objective.

\begin{landscape}

\begin{table*}[t]
\caption{BARON with bounds and cuts. PT, BT and TT represent times of preprocessing, BARON solution and total computation in seconds.  {RG represents the  percentage root gap.}}
\label{BARON cut}
\begin{center}
\begin{tabular}{ccrrrrrrrrrrrrr}
\hline
           &            &        \multicolumn{ 2}{c}{BARON}          & &            \multicolumn{ 4}{c}{BARON with bounds}                &  &           \multicolumn{5}{c}{BARON with bounds and cuts} \\
\cline{3-4} \cline{6-9} \cline{11-15}
  based on &       type &  BT (s) & {RG (\%)} & & PT (s) &  BT (s) &  TT (s)     & {RG (\%)} &   & PT (s) &  BT (s) &  TT (s) &      \#cuts  & {RG (\%)} \\
\hline
     case9 &  quadratic &  1.17 &  11.72 &  &     4.34 &       1.08 &       5.42 & 9.72 &    &   4.41 &       1.01 &       5.42 &          6 & 9.71\\

     case9 &     linear &   1.11 &16.13 &  &     4.12 &       0.86 &       4.98 &    16.02&   &    4.42 &       1.00 &       5.42 &          6& 16.05\\

    case9Q &  quadratic & 1.36 & 16.91&  &     4.36 &       1.22 &       5.58 & 8.32  &     &  4.34 &       1.11 &       5.46 &          6 &7.50\\

    case14 &  quadratic &  35.32 & 9.98& &      7.11 &      30.46 &      37.56 &  1.98 &   &    6.89 &      41.99 &      48.88 &          7 &1.98\\

    case14 &     linear &  0.79 & 0.21&   &    6.95 &       0.83 &       7.79 &  0.15&   &    6.89 &       0.91 &       7.80 &          7 &0.39\\

case\_ieee30 &  quadratic &  8347.79  &46.93 &   &   16.91 &     900.50 &     917.41 &  29.88  &  &    17.28 &       0.36 &      17.63 &         14 &0.00\\

case\_ieee30 &     linear & 2494.31 &46.67 &   &   16.89 &     249.48 &     266.37 &   33.69 &   &   16.96 &       0.34 &      17.30 &         14& 0.00\\

    case30 &  quadratic &  2.52  &9.13 &  &    17.21 &       1.91 &      19.12 &  7.74 &     & 16.94 &       4.42 &      21.35 &         13&7.39 \\

    case30 &     linear & 8.50  &5.79 &    &  17.60 &       2.40 &      19.99 &  5.15 &      &16.23 &       1.91 &      18.14 &         13&4.41 \\

   case30Q &  quadratic &  5.16  & 13.25&   &   16.53 &       2.39 &      18.93 &   12.11&     & 16.80 &       1.83 &      18.64 &         13 &4.27\\

    case39 &  quadratic & 110.59 & 8.89&   &   28.07 &      26.03 &      54.10 & 0.48  &      &27.72 &      33.25 &      60.98 &         12&1.12 \\

    case39 &     linear &  1566.88  &2.56 &    &  26.94 &      72.80 &      99.74 & 0.51  &    &  28.17 &      42.62 &      70.79 &         12 &0.52\\

  {  case57 } &  quadratic & $>10800$ & 46.69&    &  41.57 &     0.66 &     42.23&  0.00 &     &40.17&     0.80 &      40.97 &         14 &0.00 \\

  {  case57 }  &     linear &  $>10800$  & 45.17&  &   42.19 &      0.67 &      42.87 &  0.00 &     & 45.23 &      0.67 &     45.91 &         14& 0.00\\
\hline
\end{tabular}

\end{center}
\end{table*}

\end{landscape}

\noindent However, compared to the case without bounds and cuts, BARON still requires less amount of time.

{For 57-bus instances, BARON without bounds was not able to certify the global optimal solution within 3 hours time limit.   However, the strengthened variable bounds and valid inequalities enable BARON to solve these instances to global optimality within only 46 seconds.}

{As a final note, we should note that the applicability of the valid inequalities  proposed in this section is not limited to the global optimization of radial networks, they can be used in meshed networks as well. 
Moreover, precisely the same valid inequalities can be used in SOCP relaxation whereas  the transformations $c_{ij}=e_ie_j+f_if_j$ and $s_{ij}= e_if_j-e_jf_i$  enable us to obtain linear matrix inequalities to be added to  SDP relaxation. Although, for our instances,  we have not observed any lower bound improvement in SOCP/SDP relaxations by the inclusion of the valid inequalities, we obtain stronger root node relaxations in BARON. 
Let RG represent the percentage root gap calculated as $100\times (1 - z_{\text r} / z_{\text g})$, where $z_{\text r}$ and $ z_{\text g}$ are respectively the values of root node relaxation for BARON and global optimal solution. As we can see from Table \ref{BARON cut} that addition of bounds and valid inequalities strengthen the root node relaxation of BARON in general. In fact,  case\_ieee30 and case57  instances are already solved at the root node. We should note that occasionally RG of BARON with bounds and cuts is slightly worse than BARON with bounds. However, this is due to the fact that valid inequalities change the problem structure and may lead to different preprocessing  procedures carried out by the solver at the root node.}

\section{Conclusions}
\label{section:conc}

In this work, we study the impact of generation lower bounds on the performance of convex relaxations of AC OPF problems. For the fundamental two-bus one-generator model, we provide a complete characterization of all possible outcomes of the SOCP relaxation together with a detailed study of the projected feasible regions of the OPF and SOCP relaxation. 
We provide a library of radial network instances {that} demonstrate large optimality gaps for SDP and SOCP relaxations. We also propose valid inequalities for the SOCP relaxation, which prove to be useful in reducing the computation time of global solver BARON. We {remind the reader} here that SDP relaxations are 
very powerful and their importance is definite. Our work only {serves to} demonstrate the limitations of SDP relaxations and emphasizes  the importance and the need to develop efficient global methods in solving OPF problems.


\appendix

\subsection{ {A Note on BARON} }
\label{app:baron}

{
Branch-and-Reduce Optimization Navigator (BARON) is a general purpose global solver designed to solve Mixed-Integer Nonlinear Programs (MINLP). The details of the algorithm can be found in \cite{Tawarmalani}. Here, we only describe the key elements in the algorithm.
}

{
Given an MINLP, BARON first transforms the problem into a  factorable form \cite{mccormick} using  compositions of sums and products of single-variable functions. This form enables the solver to find polyhedral outer approximations of the feasible region using McCormick envelopes or other envelopes, depending on the type of the nonlinearity. After this step, linear programs are used to solve the convex relaxations.  Also, NLP solvers are utilized to find feasible solutions, which  serves as upper bounding heuristics.
}

{
There are two features in the algorithm which help to refine the convex relaxation by updating bounds on the variables. First of these features is called ``range reduction". If a variable is at its upper bound in a relaxed subproblem, its lower bound may be improved.  In case reduction step fails, then the second feature called ``spatial branching" is applied. In particular, a variable is selected and bisection is done to create two new nodes to be explored with new bounds on the branched variable.  The algorithm terminates when either all nodes are explored or lower bound and upper bound are close enough to a certain degree, which is predefined by the user.
}

{
As this brief explanation demonstrates, the success of BARON heavily depends on the quality of convex relaxations. In particular, good variable bounds are crucial in order to obtain tight polyhedral relaxations. This will help the solver to limit the number of branching steps and reduce the depth of search tree. Similarly, problem-specific valid inequalities can be useful to  tighten the convex relaxation.
}


\bibliographystyle{plain}
\bibliography{references}

\begin{thebibliography}{10}

\bibitem{Mosek}
{\em MOSEK Modeling Manual}.
\newblock MOSEK ApS, 2013.

\bibitem{bose2012quadratically}
Subhonmesh Bose, Dennice~F Gayme, K~Mani Chandy, and Steven~H Low.
\newblock Quadratically constrained quadratic programs on acyclic graphs with
  application to power flow.
\newblock {\em arXiv preprint arXiv:1203.5599}, 2012.

\bibitem{bose2011optimal}
Subhonmesh Bose, Dennice~F Gayme, Steven Low, and K~Mani Chandy.
\newblock Optimal power flow over tree networks.
\newblock In {\em {49th Annual Allerton Conference on Communication, Control,
  and Computing (Allerton)}}, pages 1342--1348, 2011.

\bibitem{Bukhsh}
Waqquas~A. Bukhsh, Andreas Grothey, Ken McKinnon, and Paul Trodden.
\newblock Local solutions of optimal power flow.
\newblock {\em IEEE Trans. on Power Syst.}, 28(4):4780 -- 4788, 2013.

\bibitem{Carpentier}
J.~Carpentier.
\newblock Contributions to the economic dispatch problem.
\newblock {\em Bulletin Society Francaise Electriciens}, 8(3):431--447, 1962.

\bibitem{Exposito99}
A.~G. Exp\'{o}sito and E.~R. Ramos.
\newblock Reliable load flow technique for radial distribution networks.
\newblock {\em IEEE Trans. on Power Syst.}, 14(3):1063 -- 1069, 1999.

\bibitem{Gan}
Lingwen Gan, Na~Li, Steven Low, and Ufuk Topcu.
\newblock Exact convex relaxation for optimal power flow in distribution
  networks.
\newblock {\em SIGMETRICS Perform. Eval. Rev.}, 41(1):351--352, June 2013.

\bibitem{Jabr06}
Rabih~A. Jabr.
\newblock Radial distribution load flow using conic programming.
\newblock {\em {IEEE} Trans. Power Syst.}, 21(3):1458--1459, 2006.

\bibitem{Jabr07}
Rabih~A. Jabr.
\newblock A conic quadratic format for the load flow equations of meshed
  networks.
\newblock {\em {IEEE} Trans. Power Syst.}, 22(4):2285--2286, 2007.

\bibitem{Jabr08}
Rabih~A. Jabr.
\newblock Optimal power flow using an extended conic quadratic formulation.
\newblock {\em {IEEE} Trans. Power Syst.}, 23(3):1000--1008, 2008.

\bibitem{josz}
C{\'e}dric Josz, Jean Maeght, Patrick Panciatici, and Jean~Charles Gilbert.
\newblock Application of the moment-sos approach to global optimization of the
  opf problem.
\newblock {\em {IEEE} Trans. on Power Syst.}, 2014.

\bibitem{lasserre}
Jean~B Lasserre.
\newblock Global optimization with polynomials and the problem of moments.
\newblock {\em {SIAM Journal on Optimization}}, 11(3):796--817, 2001.

\bibitem{Lavaei14}
J.~Lavaei, D.~Tse, and Baosen Zhang.
\newblock Geometry of power flows and optimization in distribution networks.
\newblock {\em {IEEE} Trans. Power Syst.}, 29(2):572--583, March 2014.

\bibitem{Lavaei12}
Javad Lavaei and Steven~H. Low.
\newblock Zero duality gap in optimal power flow problem.
\newblock {\em {IEEE} Trans. on Power Syst.}, 27(1):92--107, 2012.

\bibitem{Lesieutre}
Bernard~C. Lesieutre, Daniel~K. Molzhan, Alex~R. Borden, and Christopher~L.
  DeMarco.
\newblock Examining the limits of the application of semidefinite programming
  to power flow problems.
\newblock In {\em Forty-Nine Annual Allerton Conference}, pages 1492--1499,
  2011.

\bibitem{low2014convexi}
Steven~H Low.
\newblock Convex relaxation of optimal power flow, part i: Formulations and
  equivalence.
\newblock {\em IEEE Trans. Control Netw. Syst.}, 1(1):15--27, March 2014.

\bibitem{low2014convexii}
Steven~H Low.
\newblock Convex relaxation of optimal power flow, part ii: Exactness.
\newblock {\em {IEEE Trans. Control Netw. Syst.}}, 1(2):177--189, June 2014.

\bibitem{mccormick}
Garth~P McCormick.
\newblock Computability of global solutions to factorable nonconvex programs:
  Part {I} -- convex underestimating problems.
\newblock {\em Mathematical Programming}, 10(1):147--175, 1976.

\bibitem{Phan}
Dzung~T. Phan.
\newblock Lagrangian duality and branch-and-bound algorithms for optimal power
  flow.
\newblock {\em Operations Research}, 60(2):275--285, 2012.

\bibitem{sojoudi2012physics}
Somayeh Sojoudi and Javad Lavaei.
\newblock Physics of power networks makes hard optimization problems easy to
  solve.
\newblock In {\em {IEEE Power and Energy Society General Meeting}}, pages 1--8,
  2012.

\bibitem{Sojoudi}
Somayeh Sojoudi and Javad Lavaei.
\newblock On the exactness of semidefinite relaxation for nonlinear
  optimization over graphs: Part ii.
\newblock In {\em IEEE 52nd Annual Conference on Decision and Control (CDC)},
  pages 1051--1057, Dec 2013.

\bibitem{Sun}
Andy Sun and Dzung~T. Phan.
\newblock {\em Wiley Encyclopedia of Operations Research and Management
  Science}, chapter Some Optimization Models and Techniques for Electric Power
  System Short-term Operations.
\newblock John Wiley \& Sons, Inc., 2013.

\bibitem{BARON}
M.~Tawarmalani and N.~V. Sahinidis.
\newblock A polyhedral branch-and-cut approach to global optimization.
\newblock {\em Mathematical Programming}, 103(2):225--249, 2005.

\bibitem{Tawarmalani}
Mohit Tawarmalani and Nikolaos~V Sahinidis.
\newblock {\em Convexification and Global Optimization in Continuous and
  Mixed-Integer Nonlinear Programming: Theory, Algorithms, Software, and
  Applications}, volume~65.
\newblock Springer, 2002.

\bibitem{wachter}
Andreas W{\"a}chter and Lorenz~T Biegler.
\newblock On the implementation of an interior-point filter line-search
  algorithm for large-scale nonlinear programming.
\newblock {\em {Mathematical Programming}}, 106(1):25--57, 2006.

\bibitem{Zhang12}
Baosen Zhang and D.~Tse.
\newblock Geometry of feasible injection region of power networks.
\newblock In {\em Communication, Control, and Computing (Allerton), 2011 49th
  Annual Allerton Conference on}, pages 1508--1515, Sept 2011.

\bibitem{Matpower}
R.D. Zimmerman, C.E. Murillo-Sanchez, and R.J. Thomas.
\newblock {MATPOWER}: Steady-state operations, planning, and analysis tools for
  power systems research and education.
\newblock {\em {IEEE} Trans. Power Syst.}, 26(1):12--19, Feb 2011.

\end{thebibliography}

\end{document}